%version as of 25/10/99
\input amstex
\documentstyle{amsppt}
\magnification=\magstep1
\pageheight{24truecm}
\vcorrection{-2truecm}

\font\bigbf=cmbx10 scaled\magstep2
\font\biggbf=cmbx10 scaled\magstep3

\NoBlackBoxes

\define\bsr{\operatorname{\text{\rm bsr}}}
\define\cl{\operatorname{\text{\rm cl}}}
\define\clr{\operatorname{\text{\rm cr}}}
\define\End{\operatorname{\text{\rm End}}}
\define\Hom{\operatorname{\text{\rm Hom}}}
\define\wtR{\widetilde{R}}

\topmatter

\title
{\biggbf An Infinite \nolinebreak Analogue \nolinebreak of \\
\\
Rings with \nolinebreak Stable \nolinebreak Rank \nolinebreak One}\endtitle

\rightheadtext{An infinite analogue of rings with stable rank one}
\author
Pere Ara, Gert K. Pedersen\\
{\rm and} Francesc Perera
\endauthor
\address{Departament de Matem\`atiques, Universitat Aut\`onoma de Barcelona,
  08193, Bellaterra (Barcelona), Spain}\endaddress
\email{para\@mat.uab.es,
  perera\@mat.uab.es} \endemail
\address{Department of Mathematics, University of Copenhagen,
  Universitetsparken 5, DK-2100 Copenhagen \O, Denmark}\endaddress
\email{gkped\@math.ku.dk, perera\@math.ku.dk}\endemail
\date{October 1999}
\enddate
 \endkeywords 

\thanks{Research supported by the European Community (HCM Programme), the
Danish Research Council (SNF), the DGESIC (Spain) and the Comissionat
per Universitats i Recerca de la Generalitat de Catalunya.} \endthanks

\abstract{Replacing invertibility with quasi-invertibility in Bass' first
  stable range condition we discover a new class of rings, the
  $QB-$rings. These constitute a considerable enlargement of the class of rings
  with stable rank one ($B-$rings), and include examples like $\End_\Bbb F
  (V)$, the ring of endomorphisms of a vector space $V$ over some field $\Bbb
  F$, and $\Bbb B(\Bbb F)$, the ring of all row- and column-finite matrices
  over $\Bbb F$. 

  We show that the category of $QB-$rings is stable under the formation of
  corners, ideals and quotients, as well as matrices and direct limits. We
  also give necessary and sufficient conditions for an extension of $QB-$rings
  to be a $QB-$ring, and show that extensions of $B-$rings often lead to
  $QB-$rings. Specializing to the category of exchange rings
  we characterize the subset of exchange $QB-$rings as those in which every
  von Neumann regular element extends to a maximal regular element, i\.e\. a
  quasi-invertible element. Finally we show that the $C^*-$algebras that are
  $QB-$rings are exactly the extremally rich $C^*-$algebras studied by
  L.G. Brown and the second author.} 
\endabstract
\endtopmatter

\vskip0.5truecm
%\newpage

\subhead\nofrills{\bigbf 1. Introduction}\endsubhead

\bigskip

According to Bass, see e.g. \cite{\bf 6}, a unital ring $R$ will have $n$
in its {\it stable range} if whenever $(a_0,\cdots,a_n)$ is a left
unimodular row, i.e. $Ra_0+\cdots+Ra_n=R$, there is a unimodular row
$(b_1,\cdots,b_n)$ of the form $b_k=a_k+y_ka_0$ for some elements
$y_k$ in $R$, $1\le k\le n$. The smallest $n$ in the Bass range is
-- through a fortuitous error in translation from English to Russian
and back -- called the {\it stable rank} of $R$. This number, denoted by
$\bsr(R)$, is important in the non-stable $K$-theory of $R$, since the natural
morphisms $GL_n(R)/E_n(R) \to K_1(R)$, where $E_n(R)$ denotes the multiplicative
group generated by the elementary $n\times n$-matrices, will be surjective
whenever  $n\ge\bsr(R)+1$ and injective whenever $n\ge\bsr(R)+2$,
cf\.\, \cite{{\bf 6}, Theorem 2.1}.  

Specializing to rings with $\bsr(R)=1$ we see that they are
characterized by the condition that whenever $Ra+Rb=R$, then
$R(a+yb)=R$ for some $y$ in $R$. But in such a {\it $B-$ring} (to
use a terminology suggested by Vaserstein, \cite{\bf 29}) any left (or
right) invertible element is actually invertible, cf\.\, \cite{{\bf
    29}, Theorem 2.6}, so the definition of a $B-$ring can be
simplified to the demand:

\smallskip

\noindent $(*)$ \quad\;\; Given $a, b$ and $x$ in $R$ such that $xa+b=1$, 

\qquad there is a $y$ in $R$ such that $a+yb\in R^{-1}$.

\smallskip

$B-$rings have many pleasant properties, notably the cancellation
property for finitely generated projective $R$-modules, which states that
if $M$, $N$ and $P$ are $R$-modules and $P$ is projective and finitely 
generated, then $M\oplus P$ is isomorphic to $N\oplus P$ if and only if $M$ 
is isomorphic to $N$. Also it is clear that the unit can not be equivalent to 
any other idempotent, so these rings are ``finite''. Through the work of 
Rieffel, cf. \cite{\bf 28}, the Bass stable range was introduced in 
$C^*$-algebra theory, and linked with \v Cech's covering dimension; and
$C^*-$algebras having stable rank one were identified with those unital 
$C^*-$algebras $A$ for which the set $A^{-1}$ of invertible elements is dense 
in $A$. In particular, the algebra $C(X)$ of complex functions on a compact
Hausdorff space $X$ is a $B-$ring if and only if $\dim X\le 1$. 

Replacing the set of units $A^{-1}$ in a $C^*-$algebra $A$ with the union of left and
right invertible elements as in \cite{{\bf 24}}, and more generally (and
technically much more demanding) with the set $A^{-1}_q$ of quasi-invertible
elements in $A$, a theory of socalled {\it extremally rich} $C^*-$algebras has
been developed in \cite {\bf 8}, \cite {\bf 9}, \cite {\bf 10},
\cite{\bf 11} and \cite{\bf 12}. These algebras are characterized
topologically by the fact that the quasi-invertible elements form a dense set,
but one of the equivalent conditions is a version of Bass' first stable range
condition, cf\. \cite{{\bf 8}, Theorem 3.3}. This we show to be equivalent
with condition $(*)$ above, when $R^{-1}$ replaced by the set $R^{-1}_q$ of
quasi-invertible elements, see Proposition 9.1.

In the present work we systematically build a theory for rings that satisfy
the condition  

\smallskip

\noindent $(**)$ \quad Given $a, b$ and $x$ in $R$ such that $xa+b=1$, 

\qquad there is a $y$ in $R$ such that $a+yb\in R^{-1}_q$.

\smallskip

\noindent These we call $QB-$rings. Our aim, as in the $C^*-$algebra project
mentioned above, is to extend as much as possible of the theory of $B-$rings
to the much larger class of $QB-$rings. In this paper we concentrate on the
categorical properties of $QB-$rings. In subsequent papers we plan to show
that the class of $QB-$rings is the proper carrier for non-stable $K-$theory
and for an index theory, cf\. \cite{\bf 27}, including a generalized index
set, cf\. \cite {\bf 10}.

In \S 2 we define quasi-invertibility, and prove that it implies a strong
form of von Neumann regularity. In fact, every quasi-invertible element is
maximal in the natural order (by extension) on von Neumann regular elements. 
These relations are treated in some detail in order, we hope, to convince the
reader that quasi-invertibility is an important concept linking von Neumann
regularity with invertibility. We then in \S 3 define $QB-$rings by condition
$(**)$ above, and show that it is left-right symmetric. It is easy to see that
the $QB-$property passes to quotients, but in order to show that it passes to
ideals we need to reformulate the concept of quasi-invertibility to make sense
for non-unital rings. This is done in Section 4 by introducing
quasi-adversibility in exactly the same manner as Kaplansky used
adversibility instead of invertibility (although not quite using that name).  

In Section 5 we study ``skew corners'', i\.e\. subsets of the form $pRq$,
where $p$ and $q$ are idempotents in $R$. We develop a suitable notion of
$QB-$corner in order to prove in Section 6 that the class of $QB-$rings is
stable under matrix formation and under Morita equivalence.

We establish in Section 7 necessary and sufficient conditions for an extension
of $QB-$rings to be again a $QB-$ring. In particular we obtain easily
verifiable sufficient conditions, when one of the rings in an extension is a
$B-$ring. One of the motivating examples here is, of course, the algebra of
compact perturbations of Toeplitz operators, which is known to be an extension
of two $B-$rings, viz\. the algebra of compact operators and the algebra of
functions on the unit circle, but which fails to be a $B-$ring
itself. However, it is a $QB-$ring. Some algebraic analogues of the
Toeplitz algebra are also presented.

Specializing in Section 8 to exchange rings we show how the quasi-invertible
elements are precisely the maximally extended von Neumann regular
elements. This leads to a characterization of $QB-$exchange rings as exactly
those exchange rings in which every von Neumann regular element can be
extended to a maximal regular element. Finally in Section 9 we prove that a
$C^*-$algebra is a $QB-$ring if and only if it is {\it extremally rich} as
defined in \cite{\bf 8}. This means that we can use all the examples described
in that and subsequent papers, in particular we can show that our extension
results are best possible. 

\vskip1truecm
%\vskip0.5truecm

\subhead\nofrills{\bigbf 2. Quasi-Invertibility}\endsubhead

\bigskip

\definition{2.1. Definitions} In this section $R$ will denote a unital
ring. We say that two elements $x$ and $y$ in $R$ are {\it centrally
  orthogonal},  in symbols  $x\perp y$, if $xRy=0$ and $yRx=0$. Similarly, two
subsets $I$ and $J$ of $R$ are centrally orthogonal if $IRJ=0=JRI$. If $I$ and
$J$ are actually ideals in $R$, then $I\perp J$ simply means that $IJ=0=JI$. We
shall refer to this situation by saying that $I$ and $J$ are {\it orthogonal
  ideals}. Note that this does not necessarily means that $I \cap J=0$, only
that $I\cap J$ is contained in the prime radical of $R$.

We have defined our orthogonality relations above so that they will apply in
any ring. In important examples $R$ will be at least semi-prime if not
semi-primitive. In that case $xRx=0$ implies that $x=0$ for any element $x$ in
$R$. Semi-primeness of $R$ is also equivalent to the demand that $I^2=0$
implies $I=0$ for any ideal $I$ in $R$, so that for any pair of ideals $I$ and
$J$, the three conditions $IJ=0,\; JI=0$ and $I\cap J=0$ are all
equivalent. In particular, $xRy=0 \iff yRx=0$ for any pair of
elements $x,y$ in a semi-prime ring $R$.

An element $u$ in an arbitrary ring $R$ is said to be {\it quasi-invertible}
if there exist elements $a,b$ in $R$ such that
 $$
(1-ua)\perp(1-bu)\;.\tag{$*$}
 $$
The set of quasi-invertible elements in $R$ will be denoted by
$R_q^{-1}$. Similarly, the sets of left invertible, right invertible
and two-sided invertible elements will be denoted by $R_{\ell}^{-1}$,
$R_r^{-1}$ and $R^{-1}$, respectively.

If $u\in R_q^{-1}$, and if $I$ and $J$ denote the ideals of $R$ generated
by $1-ua$ and $1-bu$, respectively, then evidently
$u+I\in(R/I)_r^{-1}$ and $u+J\in(R/J)_{\ell}^{-1}$, whereas $I\perp J$ by
$(*)$.

Conversely, if $I$ and $J$ are orthogonal ideals in $R$, and $u\in R$
such that $u+I\in(R/I)_r^{-1}$ and $u+J\in(R/J)_{\ell}^{-1}$,  then
$1-ua\in I$ and $1-bu\in J$ for some elements $a,b$ in $R$, whence
$u\in R_q^{-1}$. These conditions, therefore, furnish an equivalent
definition of quasi-invertibility. In particular we see that if $R$
is a prime ring, then
$R_q^{-1}=R_{\ell}^{-1}\cup R_r^{-1}$. 

Recall that an element $a$ in a ring $R$ is 
{\it von Neumann} regular if $a=axa$ for some $x$ in $R$. We say that $x$ is a
{\it partial inverse} for $a$. Replacing if necessary $x$ with $xax$ we may
assume that $xax=x$, so that $a$ is also a partial inverse for $x$. Note that
$p=ax$ and $q=xa$ are idempotents in $R$ satisfying $pR=aR$ and
$Rq=Ra$. Conversely, if $p$ is an idempotent in $R$ such that $aR=pR$, then
$a$ is von Neumann regular, and if $p=ax$ then $x$ is a partial inverse for
$a$. 

If $u\in R_q^{-1}$, then by $(*)$ we have the equation
$(1-ua)u(1-bu)=0$. Taking $v= a+b-aub$ this implies that $u=uvu$, so
that $u$ is von Neumann regular in $R$ with partial inverse $v$. However, 
being quasi-invertible is much more than just having a partial inverse,
i.e. being von Neumann regular. For this reason (and with due apologies to
previous authors) we will in this paper reserve the name {\it quasi-inverse}
for an element  satisfying the stronger conditions in Proposition 2.2. By
computation $1-uv=(1-ua)(1-ub)$ and $1-vu=(1-au)(1-bu)$, so that we have the
relation  
 $$
(1-uv) \perp (1-vu)\;.\tag{$**$}
 $$
Moreover, replacing if necessary $v$ by $vuv$, we see that $v$ is also von
Neumann regular with partial inverse $u$. This replacement will not affect the
orthogonality relations, since $1-u(vuv)=(1-uv)(1+uv).$

We summarize our observations in the following statement:
\enddefinition

\proclaim{2.2. Proposition} Each element $u$ in $R_q^{-1}$ is von
Neumann regular, and we may choose a quasi-inverse $v$ for $u$ such that $u$
and $v$ are partial inverses for each other and the two idempotents $1-uv$
and $1-vu$ are centrally orthogonal in the sense of $(**)$. In particular,
$v\in R_q^{-1}$. \hfill $\square$  
\endproclaim 

\medskip

Note that if $u\in R_{\ell}^{-1}$ with left inverses $v$ and $v'$, then
$v'=v+v'(1-uv)$, and any element of the form $v+a(1-uv)$ will be a
left inverse for $u$. In particular, we need not have $v'=v$. We do
not, therefore, expect any unicity for the quasi-inverses of elements
in $R_q^{-1}$, so the relations described in the following result are actually
much more powerful than one might have expected.

\medskip

\proclaim{2.3. Theorem} If $u$ and $v$ are elements in $R^{-1}_q$ and
quasi-inverses for each other, so that $(1-uv)\perp (1-vu)$, then each element
of the form  
 $$
v'=v+a(1-uv)+(1-vu)b\;,\tag{$*$}
 $$
with $a,b$ in $R$, will be a quasi-inverse for $u$ in $R_q^{-1}$ and
satisfy the relations
 $$
(1-uv')\perp (1-v'u)\,,\quad (1-uv')\perp (1-vu)\,,\quad (1-uv)\perp (1-v'u)\;.
 $$
 
Conversely, if $v'$ is any partial inverse for $u$, then $v'\in R_q^{-1}$
and has the form $(*)$ with $a=b=v'$. Moreover,
 $$
1-uv=(1-uv)(1-uv')\quad\text{and}\quad 1-uv'=(1-uv')(1-uv)\;,\tag{$**$}
 $$
and similarly for the $vu$ and $v'u$ products. In particular, the idempotents
$1-uv$ and $1-uv'$ are Murray-von Neumann equivalent, as are $1-vu$ and
$1-v'u$, so there are orthogonal ideals $I$ and $J$ in $R$, such that
$uv=1=uv'$ modulo $I$ and $vu=1=v'u$ modulo $J$. 
\endproclaim

\demo{Proof} Evidently each $v'$ in $R$ of the prescribed form
$(*)$ will be a partial inverse for $u$. Moreover,
 $$
1-uv'= (1-ua)(1-uv)\quad\text{and}\quad 1-v'u = (1-vu)(1-bu)\;,
 $$
 whence $(1-uv')\perp (1-v'u)$, so that $v'\in R_q^{-1}$ and $v'$ is a
 quasi-inverse for $u$.

Conversely, if $uv'u=u$ for some $v'$ in $R$ then
 $$
0=(1-vu)v'(1-uv)=v'-vuv'-v'uv+v\;.
 $$
 Consequently,
 $$
v'=v+2v'-vuv'-v'uv=v+(1-vu)v'+v'(1-uv)\;,
 $$
 as desired, and $v'\in R_q^{-1}$ with $(1-uv')\perp (1-v'u)$ by the
first part of the proof.

The equations in $(**)$ follow by straightforward computations.
The equivalence between the idempotents $1-uv$ and $1-uv'$ is evident from  
$(**)$, and therefore $1-uv$ and $1-uv'$ generate the same ideal $I$.
Similarly, $1-vu$ and $1-v'u$ generate the same ideal $J$. These
ideals are orthogonal since $(1-uv)\perp (1-vu)$, and evidently the
desired relations are satisfied in $R/I$ and $R/J$, respectively.
\hfill$\square$
\enddemo

\medskip

\definition{2.4. Definitions} If $a$ and $b$ are von Neumann regular elements
in $R$ we say that $b$ {\it extends} $a$, and write $a\le b$, if
 $$
a=axb=bxa=axa\tag{$*$}
 $$
for some $x$ in $A$. Thus $x$ is a partial inverse for $a$, and we may
assume that $a$ is also a partial inverse for $x$.

Taking $p=ax$ and $q=xa$ in the equations above we see from $(*)$
that $p$ and $q$ are idempotents in $R$ such that $pb=bq=a$.
Moreover, $pR=aR$ and $Rq=Ra$, since $x$ is a partial inverse for $a$.
Conversely, if $a\in R$ and if we can find idempotents $p$ and $q$, and a von
Neumann regular element $b$, such that 
$$
pR=aR\quad\text{and}\quad Rq=Ra,\quad \text{and moreover}\quad pb=bq=a,\tag{$**$}
$$ 
then $a$ is von Neumann regular and $a\le b$. Thus, $(**)$ furnishes an
alternative description of the relation $\le$. Observe from this that
if $a\le b$ and $aR=bR$ (or if $Ra=Rb$), then $a=b$ (because then $pR=bR$, so
$b=pb=a$).   

The idempotents $p$ and $q$, above, depend not only on $a$ but also
on $x$. Nevertheless we may think of them as the ``range'' and the
``source'' of $a$, noting that $a$, as a left multiplier on $R$, is a
bijection from $qR$ onto $pR$ with kernel $(1-q)R$ and cokernel
$(1-p)R$. In this setting $a\le b$ expresses an ordinary
extension of operators. 

If $y$ is a partial inverse for $b$, with $b$ as its partial inverse,
then from $(*)$ we derive the simple relations
\roster
\item"(i)" $ax=bx$, $xa=xb$.
\item"(ii)" $a=ayb=bya$.
\item"(iii)" $a=bxb=aya$.
\endroster 

Now put $x'=yay$ and check that
 $$
ax'a=a,\quad x'ax'=x',\quad x'=x'ay=yax'.
 $$
 Thus, replacing if necessary $x$ with $x'$, we may assume that also
$x\le y$.

This implies that $\le$ is a transitive relation. For if
$a\le b$ and $b\le c$ we may assume that $a$ and $b$ have
partial inverses $x$ and $y$, respectively, such that $x\le y$ as above.
Then
 $$
cxa=c(yax)a=(cy)(axa)=(by)a=a,
 $$
 and similarly $axc=a$, so $a\le c$.

If $a\le b$ and $b\le a$ then in particular $aR=bR$, whence $a=b$.
 We summarize our observations in the following result:
\enddefinition

\medskip

\proclaim{2.5. Proposition} The relation $\le$ defined above on
the set $R^r$ of von Neumann regular elements in a ring $R$ gives a
partial order on $R^r$. \hfill$\square$
\endproclaim

\medskip

\proclaim{2.6. Lemma} Let $R$ be a unital ring, and let $a$ be a
regular element in $R$. Write $aR=pR$ and $Ra=Rq$ for some idempotents 
$p$ and $q$. If $b$ is a regular element such that $pb=a=bq$, and if
$q'$ is an idempotent satisfying $Rq=Rq'$, then there is an element $b'$ 
such that $Rb=Rb'$, $bR=b'R$, and $pb'=a=b'q'$.
\endproclaim

\demo{Proof} Let $b'=a+(1-p)b(1-q')$. Clearly then, $pb'=pa=a$ and
$b'q'=aq'=a$. Note that we can also write $b'=a+b(1-q')-pb(1-q')=
a+b(1-q')=bq+b(1-q')$. Since $qq'=q$ and $q'q=q'$, it follows that
$w=q+(1-q')$ is invertible with $w^{-1}=q'+(1-q)$. From above we have that
$b'=bw$, whence $b=b'w^{-1}$, and thus $bR=b'R$. 

Since $Rq'=Ra\subset Rb$, we have that $b'=a+(1-p)b-(1-p)bq'\in Rb$,
whence $Rb'\subset Rb$. Write $q'=ta$, with $q'tp=t$, and compute that
$$
\aligned
(bt+(1-p))b' &=(bt+(1-p))(a+(1-p)b(1-q'))\\
             &=bta+(1-p)b(1-q')=bq'+b(1-q')=b\;,
\endaligned
$$
so $Rb\subset Rb'$, hence $Rb=Rb'$. \hfill $\square$
\enddemo

\medskip

\proclaim{2.7. Lemma} If $a\le b$ in $R^r$ there exist idempotents $p$
and $q$ in $R$ such that $pb=a=bq$ and $b-a\in (1-p)R(1-q)$. Conversely, if
$a$ and $c$ are elements in $R^r$, where $a\in pRq$ and $c\in(1-p)R(1-q)$
for some idempotents  $p$ and $q$ in $R$, then $a\le a+c$ in $R^r$.
\endproclaim

\demo{Proof} By $(**)$ in 2.4 we have idempotents $p$ and $q$ such that
$pb=a=bq$ and $(1-p)a=0=a(1-q)$. Consequently,

$b=(p+1-p)b(q+1-q)=a+(1-p)a+a(1-q)+(1-p)b(1-q)=a+(1-p)b(1-q)$. 

Conversely, if $a$ and $c$ are regular elements in $pRq$ and $(1-p)R(1-q)$,
respectively, we can find partial inverses $x$ and $z$ for them. We may assume
that $x\in qRp$ and $z\in (1-q)R(1-p)$, replacing them if necessary with
$qxp$ and $(1-q)z(1-p)$. Then with $b=a+c$ we have $bxa=a=axb$, so $a\le b$. 
Moreover, $y=x+z$ will be a partial inverse for $b$ (even satisfying $x\le
y$). \hfill $\square$ \enddemo

\medskip

\proclaim{2.8. Proposition} For a unital ring $R$, each element in  $R^{-1}_q$
is maximal in $R^r$ with respect to the ordering $\le$
\endproclaim

\demo{ Proof} If $u\in R_q^{-1}$ and $u\le a$ for some $a$ in
$R^r$, then $u=uvu=avu=uva$ for some $v$ in $R$. But then $v$ is a
quasi-inverse for $u$ by Theorem 2.3. In particular,
$$
a= a-(1-uv)a(1-vu) =uva+avu-uvavu= u+u-uvu=u\;. 
$$
\hfill$\square$
\enddemo

\proclaim{2.9. Proposition} For a von Neumann regular element  $a$ in a
unital ring $R$ the following conditions are equivalent:
\roster
\item"(i)"   $a\le u$ for some $u$ in $R^{-1}_q$\,;
\item"(ii)"  $a=ava$ for some $v$ in $R^{-1}_q$\,.
\endroster
\endproclaim

\demo{ Proof} (i) $\implies$ (ii) If $a\le u$ in $R^r$, choose a quasi-inverse
$v$ for $u$. By (iii) in 2.4 this implies that $a=ava$. Indeed, $a$ has a
partial inverse $x$ such that $a=axu=uxa$, whence
$$
a=a(xa)=(axu)(xa)=(ax)(uvu)(xa)=(axu)v(uxa)=ava\,.
$$

\noindent (ii)$\implies$ (i) If $a=ava$ with $v$ in $R^{-1}_q$, consider the
idempotents $p=va$ and $q=av$. Let $w$ be a partial inverse for $v$ and
consider the idempotents $e=vw$ and $f=wv$. Then $e'=(1-p)e$ is an idempotent
because $ep=vwva=va=p$. Evidently $e'\le 1-p$. Moreover,
$(p+e')v=pv+ev-pev=v$, so $vR\subset (p+e')R = pR + e'R$. On the other hand, both
$pR\subset vR$ and $e'R \subset vR$, so we have equality. By symmetry we can
find an idempotent $f' \le 1-q$, such that $v(q+f')=v$ and $Rv=Rq+Rf'$. Write
$p+e'=vs$ and $q+f'=tv$ for some elements $s, t$ in $R$. Then put $u=tvs$ and
check that $vu=vtvs=v(q+f')s=vs=p+e'$ and similarly $uv=q+f'$. In particular,
$vuv=v$, so $u\in R^{-1}_q$ by Theorem 2.3. Finally, $avu=a(p+e')=ava+0=a$
and $uva=(q+f')a=ava+0=a$, so that $a\le u$, as desired. \hfill$\square$
\enddemo

\medskip

\example{2.10. Remark} Notice the affinities between the statements in Theorem
2.3 and Proposition 2.9: If $uvu=u$ for some $u$ in $R^{-1}_q$, then
necessarily $v\in R^{-1}_q$. But if $ava=a$ for some $v$ in $R^{-1}_q$, then
at least $a$ extends to some $u$ in $R^{-1}_q$.
\endexample

\vskip1truecm
%\vskip0.5truecm
%\newpage

\subhead\nofrills{\bigbf 3. QB-Rings}\endsubhead

\bigskip

\definition{3.1. Definition} For each subset $A$ of a unital ring $R$ we
define $\cl(A)$ to be the set of elements $a$ in $R$, 
such that whenever $xa+b=1$ for some elements $x$ and $b$ in $R$,
there is an element $y$ in $R$ such that $a+yb\in A$. Equivalently, $a\in
\cl(A)$ if $(a+Rb)\cap A\ne \emptyset$ whenever $Ra+Rb=R$. 
\enddefinition

\medskip

\proclaim {3.2. Lemma} The operation $\cl$ defined in 3.1 has the following
properties relative to any subsets $A$ and $B$ of $R$:
\roster
\item "{\bf (i)}" $\cl(\emptyset)=\emptyset$ and $\cl(R)=R$\,;
\item"{\bf (ii)}" $A\subset B$ implies $\cl(A)\subset \cl (B)$\,;
\item"{\bf (iii)}" $A\subset \cl (A) = \cl(\cl(A))$\,;
\item"{\bf(iv)}" If $A\ne\emptyset$ then $\Cal J(R)\subset\cl(A)$, where $\Cal J(R)$
  is the Jacobson radical of $R$\,; 
\item"{\bf (v)}" $\Cal J(R)=\cl(0)$\,;
\item"{\bf (vi)}" $\cl(A)\cap R^{-1}_\ell=A\cap R^{-1}_\ell$\,;
\item"{\bf(vii)}" $B\cl(A)\subset\cl(BA)$\,;
\item"{\bf (viii)}" If $B\subset R^{-1}$ then $\cl(A)B\subset\cl(AB)$\,;
\item"{\bf (ix)}" $\cl(A)+B\subset \cl(A+RB)$\,;
\item"{\bf (x)}" If $RB\subset B$ and $A+B\subset\cl(A)$, then
  $\cl(A)+B\subset\cl(A)$\,;
\item"{\bf (xi)}" If $\pi\colon R\to S$ is any quotient morphism, then
  $\pi(\cl(A))\subset \cl(\pi(A))$\,.
\endroster
\endproclaim

\demo{Proof} Properties {\bf (i)} and {\bf (ii)} are trivial to verify, as well
as the first part of {\bf (iii)}. To complete the argument take any $a$ in
$\cl(\cl(A))$, and assume that $xa+b=1$. Then by assumption $a+yb\in\cl(A)$
for some $y$ in $R$. But since $x(a+yb)+(1-xy)b=1$ this implies that
$a+yb+z(1-xy)b\in A$ for some $z$ in $R$. Consequently $a+(y+z-zxy)b\in A$,
whence $a\in\cl(A)$. 

\noindent{\bf(iv)} If $A\ne\emptyset$ and $z\in \Cal J(R)$, then any equation
$xz+b=1$ implies that $1-b\in \Cal J(R)$, whence $b\in R^{-1}$. Evidently then,
$A\subset z+Rb$, so $z\in \cl(A)$\,.

\noindent{\bf (v)} From (iv) we know that $\Cal J(R)\subset \cl(0)$. But if
$a\notin \Cal J(R)$ then $a\notin L$ for some maximal left ideal $L$ of $R$. As
$Ra+L=R$ by maximality, we have $xa+l=1$ for some $l$ in $L$. If $a\in \cl(0)$
this would imply that $a+yl=0$ for some $y$ in $R$, whence $a\in R$, a
contradiction. Therefore $(R\setminus \Cal J(R))\cap \cl(0)=\emptyset$, so
$\cl(0)=\Cal J(R)$\,.

\noindent{\bf (vi)} If $a\in \cl(A)\cap R^{-1}_\ell$ we may consider the
trivial decomposition $Ra+R0=R$. By assumption $a+y0\in A$, i\.e\. $a\in A$\,.

\noindent{\bf (vii)} If $e\in B$ and $a\in \cl (A)$, consider any equation
$xea+b=1$. Then $a+yb\in A$ for some $y$, whence $ea+eyb\in eA\subset BA$, and
so $ea\in\cl(BA)$\,. 

\noindent{\bf (viii)} Consider now an equation $xae+b=1$. Since $B\subset
R^{-1}$ this rewrites as $exa + ebe^{-1} =1$. Therefore $a+yebe^{-1}\in A$ for
some $y$ in $R$, whence $ae+yeb\in Ae\subset AB$, showing that $ae\in\cl(AB)$\,.

\noindent{\bf (ix)} Suppose that $e\in B$ and $a\in \cl(A)$. For any equation
$x(e+a)+b=1$ we can then find $y$ such that $a+y(xe+b)\in A$. But then 
$e+a+yb\in\ A +(1-yx)e\subset A+RB$, whence $e+a\in\cl(A+RB)$\,.

\noindent{\bf (x)} By (ix), (ii) and (iii) we have 
$$
\cl(A)+B\subset \cl(A+RB)\subset\cl(A+B)\subset\cl(\cl(A))=\cl(A)\,.
$$
             
\noindent{\bf (xi)} If $a\in \cl(A)$ and $\pi(x)\pi(a)+\pi(b)=\pi(1)$ in $S$,
then $xa+b+t=1$ in $R$ for some $t$ in $\ker\pi$. But then $a+y(b+t)\in A$ for
some $y$ in $R$, whence $\pi(a)+\pi(y)\pi(b)\in\pi(A)$, proving that
$\pi(a)\in \cl(\pi(A))$\,. 
\hfill $\square$ \enddemo

\medskip

\example{3.3. Remark} If $R$ is commutative and $Ra+Rb_i=R$ for $i= 1,2$, then
also $Ra+Rb_1b_2=R$. Moreover,
$$
a+Rb_1b_2\subset (a+Rb_1)\cap (a+Rb_2)\,.
$$
This means that the sets of ``neighbourhoods'' of $a$, each of the form
$\Cal O_a(b)=a+Rb$ for some $b$ such that $Ra+Rb=R$, is directed by
inclusion. Observe also that if $c=a+yb\in\Cal O_a(b)$, then $Rc+Rb=R$ and 
$$
\Cal O_c(b)=c+Rb=a+Rb=\Cal O_a(b)\,.
$$
Thus the sets $\Cal O_a(b)$ form the neighbourhood basis in a topology on $R$
for which $\cl$ is the closure operation. In particular, $\cl(A\cup B)=
\cl(A)\cup \cl(B)$.

For non-commutative rings this fails already when $R=\Bbb M_2(\Bbb
R)$. Nevertheless the operation $\cl$ may with advantage be compared to  a
closure. A striking case occurs in $C^*-$algebra theory, cf\. Proposition 9.1.
\endexample

\medskip

\definition{3.4. Definition} We shall be (almost) exclusively concerned with
applying the operation $\cl$ to the set $R^{-1}_q$ of quasi-invertible elements
in a unital ring $R$. Since $R^{-1}R^{-1}_q=R^{-1}_q$ and
$R^{-1}_qR^{-1}=R^{-1}_q$ we see from (vii) and (viii) in Lemma 3.2 that
$R^{-1}\cl(R^{-1}_q)=\cl(R^{-1}_q)$ and $\cl(R^{-1}_q)R^{-1} =\cl(R^{-1}_q)$.\newline 
\newline
\centerline{If $\cl(R_q^{-1})=R$ we say that $R$ is a {\it $QB-$ring}.} \newline
\newline
As mentioned in 3.1 the condition that $a\in\cl(R_q^{-1})$ is
equivalent to the demand that whenever $(a,b)$ is a left unimodular row,
i.e. $Ra+Rb=R$, then $a+yb\in R_q^{-1}$ for some $y$ in $R$. Replacing
quasi-invertibility by honest invertibility we are back at the
definition of Bass stable rank $1$, which in this setting says that
$\cl(R^{-1})=R$. Thus, $QB-$rings are a generalization of $B-$rings, and
actually a substantial weakening of this concept in the non-commutative
case. Of course, if $R$ is commutative, $R^{-1}_q=R^{-1}$. 

Evidently the definition of $\cl(R_q^{-1})$ is not left-right symmetric, so we
define $\text{\rm cr}(R_q^{-1})$  to be the set of elements $a$ in $R$, such
that whenever  $ax+b\in R^{-1}$ for some elements $x,b$ in $R$, then $a+by\in
R_q^{-1}$ for some $y$ in $R$. We have no reason to believe that
$\cl(R_q^{-1})=\text{\rm cr}(R_q^{-1})$ in general. (But see Corollary 9.2.)
However, just as for $B-$rings we do have complete symmetry in the definition
of a $QB-$ring. This is proved by adapting \cite{{\bf 13}, Lemma 1} to our
present situation. 
\enddefinition

\medskip

 \proclaim{3.5. Lemma} If $a\in R$ such that $ax+b=1$ for some $x$ in
$\cl(R_q^{-1})$ and $b$ in $R$, then $a+by\in R_q^{-1}$ for some $y$
in $R$.
 \endproclaim

\demo{Proof} Since $x\in\cl(R_q^{-1})$ we have $x+cb\in R_q^{-1}$
for some $c$ in $R$. Thus
 $$
(1-(x+cb)z)\perp (1-z(x+cb))
 $$
 for some $z$ in $R$.

Define $y=z(1-ca)$ and $d=x+(1-xa)c$. By straightforward, albeit
lengthy computations, using that $ax+b=1$, we then see that
 $$
\aligned
1-(a+by)d &= 1-(a+bz(1-ca))(x+(1-xa)c)\\
 = \cdots &= b(1-z(x+cb))(1-ac)\;.
\endaligned
 $$
 Similarly,
 $$
\aligned
1-d(a+by) &= 1-(x+(1-xa)c)(a+bz(1-ca))\\
          &= (1-xa)(1-(x+cb)z)(1-ca)\;.
\endaligned
 $$
 It follows that
 $$
(1-(a+by)d)\perp (1-d(a+by))\;,
 $$
 which shows that $a+by\in R_q^{-1}$, as desired. \hfill$\square$
\enddemo

\medskip

\proclaim{3.6. Theorem} In any unital ring $R$ we have
$$
\cl(R_q^{-1})=R \quad \iff \quad \text{\rm cr}(R_q^{-1})=R\,.
$$ 
\endproclaim

\demo{Proof} By symmetry it suffices to show that $\text{\rm cr}(R_q^{-1})=R$,
assuming that $\cl (R^{-1}_q)=R$. But that is immediate from Lemma 3.5.
\hfill $\square$ 
\enddemo

\medskip

\proclaim{3.7. Proposition}  If $I$ is an ideal in a unital ring $R$ and 
$\pi\colon R\to R/I$ denotes the quotient morphism, then
$$
\pi(R_q^{-1})\subset(R/I)_q^{-1}\quad\text{and}\quad
\pi(\cl(R_q^{-1}))\subset\cl((R/I)_q^{-1})\;.
$$
\endproclaim

\demo{Proof} Since central orthogonality is preserved under quotient morphisms
it is evident that $\pi(R_q^{-1})\subset(R/I)_q^{-1}$. The other inclusion
follows from Conditions (xi) and (ii) in Lemma 3.2. \hfill$\square$
\enddemo

\medskip

\proclaim{3.8. Corollary} If $I$ is an ideal in a unital
ring $R$ such that $I+\cl(R_q^{-1})=R$ then $R/I$ is a $QB-$ring. In
particular, the property of being a $QB-$ring passes to quotients.
\hfill $\square$ \endproclaim

\medskip

The next result, a generalization of \cite{{\bf 11}, Proposition 2.6},
indicates the position of the $B-$rings as the ``finite'' members in the
category of $QB-$rings.

\medskip
 
\proclaim{3.9. Proposition} Let $R$ be a unital $QB-$ring. Then 
the following conditions are equivalent:
\roster
\item "(i)"  $R$ is a $B-$ring \,,
\item "(ii)" $R^{-1}_q=R^{-1}$\,.
\endroster
\endproclaim

\demo{Proof} (i)$\implies$ (ii) Consider $u$ in $R^{-1}_q$ with quasi-inverse
$v$, and let $I$ and $J$ denote the orthogonal ideals generated by
$1-uv$ and $1-vu$, respectively. Then $u$ is right invertible in $R/I $, and
since  $\bsr (R/I)=1$ because $\bsr (R) =1$, this means that $u$ is actually
invertible in $R/I$ with (the image of) $v$ as its inverse. It follows that
$1-uv\in I$ and $1-vu\in I$. Similarly $1-uv\in J$ and $1-vu\in J$. Since both
elements are idempotents and $IJ=0$, the elements must be zero, and we
conclude that $uv=vu=1$, so that $u\in R^{-1}$.  

\noindent (ii)$\implies$(i) If $R$ is a $QB-$ring, then $\cl (R^{-1}_q)=R$. If in
addition $R^{-1}_q=R^{-1}$ then by definition $R$ is a $B-$ring. \hfill$\square$
\enddemo 

\medskip

Recall that a simple, unital ring $R$ is said to be {\it purely infinite} if
$R$ is not a division ring, but for any non-zero element $x$ in $R$ there are
$s, t$ in $R$ such that $sxt=1$. It is an open problem whether such a ring
must be an exchange ring, cf\. \S 8, but it is certainly well supplied
with idempotents.

\proclaim{3.10. Proposition} Every simple, purely infinite ring $R$ is a
$QB-$ring, but has infinite Bass stable rank.
\endproclaim

\demo{Proof} Take $a$ in $R$ and assume that $xa+b=1$ for some elements $x,b$
in $R$. If $b= 0$ then $a$ is left invertible, so $a \in R^{-1}_q$. Otherwise,
since $R$ is purely infinite, we can write $sbt=1$ and compute  
$$
(a+(1-at)sb)t= at+(1-at)=1\,.
$$ 
With $y=(1-at)s$ this proves that $a+yb$ is right invertible. Thus 
$a\in \cl(R^{-1}_q)$, and $R$ is a $QB-$ring.

Since $R$ is not a division ring, there is a non-zero element $z$ in
$R$ which is not left invertible. Write $czd=1$ for some elements
$c,d$ in $R$. Let $e=zdc$, and note that $e$ is a non-trivial
idempotent in $R$, equivalent to $1$. Set $f_1=1-e$. Take any non-zero
element $z_1$ in $f_1Rf_1$, which is neither right nor left
invertible. Since $R$ is purely infinite, there are elements $c_1$ and 
$d_1$ in $R$ such that $c_1z_1d_1=1$. We may assume that $c_1\in
Rf_1$. Let $f=z_1d_1c_1$. Hence $e$ and $f$ are two orthogonal
idempotents in $R$, both equivalent to $1$. Now an argument similar to
\cite{{\bf 28}, Proposition 6.5} shows that $R$ has infinite stable rank.
\hfill $\square$
\enddemo

\medskip

\example{3.11. Example} Let $V$ be a vector space over a field $\Bbb F$
and consider the ring $R=\End_\Bbb F(V)$ of all endomorphisms of $V$,
which is known to be a unital, prime ring. Then $R$ is a $QB-$ring (whereas the
Bass stable rank is $\infty$ if $V$ is infinite dimensional).

To show this, consider $a,b$ and $x$ in $R$ such that $xa+b=1$.
Evidently, then,  $\ker a\cap\ker b=0$, so we can choose a subspace
$U$ of $V$ such that $V=\ker a\oplus\ker b\oplus U$. We can also
choose subspaces $W$ and $W'$ of $V$ such that $V=b(V)\oplus
W=a(V)\oplus W'$. Define $c$ in $R$ by $c\,|W=0$ and $c|b(V)=(b\,|\ker
a\oplus U)^{-1}$. We can also define $d$ in $R$ by choosing an
injective morphism $z:\ker a\to W'$ if $\dim(\ker a)\le\dim W'$ or a
surjective morphism $z:\ker a\to W'$ if $\dim(\ker a)\ge\dim W'$, and
then take $d\,|\ker b\oplus U=0$ and $d\,|\ker a=z$. Finally, let
$y=dc$. Then $a+yb$ gives a bijective morphism from
$\ker b\,\oplus\, U$ onto $a(V)$, whereas $a+yb\,|\ker a=dcb\,|\ker
a=z$. It follows that $a+yb$ is either injective or surjective on
$V$, whence $a+yb\in R_{\ell}^{-1}\cup R_r^{-1}=R_q^{-1}$, since $R$
is prime. Consequently $a\in\cl(R_q^{-1})$, whence $R=\cl(R_q^{-1})$,
so $R$ is a $QB-$ring.
\endexample

\medskip

In later sections, notably \S 7, \S 8 and \S 9, we shall present many more
examples of $QB-$rings, as well as some counterexamples. For
the time being we just notice that when trying to generalize the above example
to the ring $\End(G)$ of an abelian group $G$ we see that if $G$ contains
$\Bbb Z$ as a direct summand, then $\End(G)$ has an idempotent $p$ such that
$p\End(G)p=\Bbb Z$. Using Corollary 5.8 it follows that
$\End(G)$ is not a $QB-$ring. 

\vskip1truecm
%\vskip0.5truecm
%\newpage

\subhead\nofrills{\bigbf 4. Rings Without Unit}\endsubhead

\bigskip

\subhead{4.1. Adversibility}\endsubhead Let $R$ be a not necessarily
unital ring. There is then a canonical unitization $\wtR=R\oplus\Bbb Z$, with
the obvious multiplication
$$
(a,n)(b,m)=(ab+ma+nb,nm)
$$
for $a,b$ in $R$ and $n,m$ in $\Bbb Z$. Evidently $\wtR$ contains
$R$ as an ideal with $\wtR/R=\Bbb Z$.

The fact that $\Bbb Z$ is not a $B-$ring and thus, being
commutative, neither a $QB-$ring, coupled with the fact that the
$QB-$property passes to quotients (Corollary 3.8) shows
that $\wtR$ will never be  a $QB-$ring. We shall therefore seek a definition
for the $QB-$property which is intrinsic for the non-unital case, even though
it may seemingly borrow some structure from $\wtR$ or other unitizations.

\medskip

Following ideas going back to Kaplansky we say that an
element $x$ in $R$ is {\it left} (respectively {\it right}) 
{\it adversible} if $a+x=ax$ (respectively $x+a=xa$) for some $a$ in $R$,  
which we call a {\it left} (respectively {\it right}) {\it adverse} to $x$. If
$x$ has both a left adverse $a$ and a right adverse $b$ it is called {\it
adversible}, in which case
$$
\align
a & = a+(x+b-xb)-a(x+b-xb)\\
& = b+(x+a-ax)-(x+a-ax)b=b\;.
\endalign
$$
The subsets of left, right and two-sided adversible elements in $R$
are denoted $R_{\ell}^\circ$, $R_r^\circ$ and $R^\circ$, respectively.

We say that an element $x$ in $R$ is {\it quasi-adversible}, in
symbols $x\in R_q^\circ$, if there exist elements $b$ and $c$ in $R$
such that
 $$
(x+b-xb)\perp (x+c-cx)\;.\tag{$*$}
 $$
Here $s\perp t$ in $R$ means that $s\wtR t=0$ and $t\wtR s=0$. We say in this
case that $s$ and $t$ are {\it centrally orthogonal}. If therefore $I$ and $J$
denote the ideals in $R$ generated by $x+b-xb$ and $x+c-cx$, respectively,
then $I\perp J$ in the sense that $s\perp t$ for any pair of elements $s$ in
$I$ and $t$ in $J$. Moreover, $x+I$ is right adversible in $R/I$ and $x+J$ is
left adversible in $R/J$. Conversely, if we can find an orthogonal pair $I$,
$J$ of ideals in $R$ (i\.e\. $I\wtR J=0$ and $J\wtR I=0$) such that
$x+I\in(R/I)_r^\circ$ and $x+J\in(R/J)_{\ell}^\circ$, then $(*)$ is satisfied,
so that $x\in R_q^\circ$. 

\proclaim{4.2. Proposition} Let $S$ be any unital ring containing
$R$ as a subring, and let $R_1$ denote the unital subring of $S$
generated by $S$ and $1$. Then we have the equalities
$$
\aligned
& 1-R_{\ell}^\circ=(R_1)_{\ell}^{-1}\cap(1-R)\,,\quad
1-R_r^\circ=(R_1)_r^{-1}\cap(1-R)\,,\\
& 1-R^\circ=(R_1)^{-1}\cap (1-R)\,,\quad1-R_q^\circ = (R_1)_q^{-1}\cap(1-R)\,.
\endaligned
$$
\endproclaim
 
\demo{Proof} Straightforward computations based on the equation\newline
\centerline{$1-(1-x)(1-a)=x+a-xa$\,.} \hfill$\square$
\enddemo

\proclaim{4.3. Corollary} If $x\in R_q^\circ$ we can find a single
element $y$ in $R$, called the quasi-adverse for $x$, such that $x+y-xy$ and
$x+y-yx$ are centrally orthogonal idempotents in $R$.
\endproclaim

\demo{Proof} Combine Propositions 2.2 and 4.2.\hfill$\square$
\enddemo

\medskip

\definition{4.4. Definition} For any ring $R$ we define
$\cl^\circ(R_q^\circ)$ to be the set of elements $a$ in $R$, such that
whenever $xa-x-a+b=0$ for some $x,b$ in $R$, there is an element $y$
in $R$ such that $a-yb\in R_q^\circ$. Note that
$R_q^\circ\subset\cl^\circ(R_q^\circ)$. If $\cl^\circ(R_q^\circ)=R$ we
say that $R$ is a $QB-$ring. Similarly, one defines
$\clr^\circ(R_q^\circ)$. This prompts the question whether the
notion of $QB-$ring in the non-unital case is also left-right
symmetric. We shall address this problem in Remark 4.8.
\enddefinition

\medskip

\example{4.5. Remarks} If $R$ is unital then
$R_q^{-1}=1-R_q^\circ$ by Proposition 4.2. Moreover, $(1-x)(1-a)+b=1$
if and only if $xa-x-a+b=0$. It follows that
$\cl(R_q^{-1})=1-\cl^\circ(R_q^\circ)$. Thus, $\cl(R_q^{-1})=R$ if and only if
$\cl^\circ(R_q^\circ)=R$, so that our two definitions for $QB-$rings coincide.

In the non-unital case $\cl^\circ(R_q^\circ)=R$ means that
$1-R\subset\cl(\wtR_q^{-1})$, which is not at all the same as 
$\wtR=\cl(\wtR_q^{-1})$.

If $R$ is an algebra over a field $\Bbb F$, we redefine $\wtR=R\oplus\Bbb F$, 
which is the minimal unitization of $R$ as an
algebra over $\Bbb F$. Now the discrepancies above vanish, and we see
that $R$ is a $QB-$algebra if and only if $\wtR$ is a $QB-$algebra.

To utilize the new definition in the non-unital case we shall
need the following reformulation of \cite{{\bf 29}, Lemma 3.5}. 
\endexample

\medskip

 \proclaim{4.6. Lemma} Let $I$ be a right ideal in a unital ring $R$.
Then for any pair of elements $a$ in $I$ and $b$ in $R$ the following
conditions are equivalent:
\roster
 \item"(i)" $1\in R(1-a)+Rb$;
 \item"(ii)" $I=I(1-a)+Ib$;
 \item"(iii)" $1-I \subset (1-I)(1-a)+Ib$;
 \item"(iv)" $xa-x-a+yb=0$ for some $x,y$ in $I$.
 \endroster
\endproclaim

\demo{Proof} (i) $\implies$ (ii) If $1=c(1-a)+db$ for some
$c,d$ in $R$ and if $t\in I$, then $t=tc(1-a)+tdb$, and $tc$, $td$
belong to $I$.

\noindent (ii) $\implies$ (iii) Given $t$ in $I$ choose $x,y$ in $I$ such
that $a-t=-x(1-a)+yb$. Then
 $$
1-t=1-a-x(1-a)+yb=(1-x)(1-a)+yb\;.
 $$

\noindent (iii) $\implies$ (iv) Choose $x,y$ in $I$ such that
$(1-x)(1-a)+yb=1$, whence $xa-x-a+yb=0$.

\noindent (iv) $\implies$ (i) If $xa-x-a+yb=0$, then $1=(1-x)(1-a)+yb$.
\hfill$\square$
\enddemo

\medskip

 \proclaim{4.7. Proposition} If $I$ is an ideal in a unital ring $R$, 
and $t\in I$, then $t\in I_q^\circ$ if and only if
$1-t\in R_q^{-1}$. Moreover, $t\in\cl^\circ(I_q^\circ)$ if and only if
$1-t\in\cl(R_q^{-1})$. 
 \endproclaim

\demo{Proof} If $t\in I_q^\circ$, then by Corollary 4.3 we can find
centrally orthogonal idempotents in $I$ of the form $p=s+t-ts$ and
$q=s+t-st$ for some $s$ in $I$. Since $I$ is an ideal in $R$ this implies that
 $$
pRq = p^2Rq^2 \subset pIq=0\;.
 $$
Similarly  $qRp=0$. It follows that $1-t\in
R_q^{-1}$ with partial inverse $1-s$.

Conversely, if $1-t\in R_q^{-1}$ we can find a quasi-inverse in $R$,
written in the form $1-s$, such that
 $$
(1-(1-t)(1-s))\perp (1-(1-s)(1-t))\;.
 $$
 The equation $1-t=(1-t)(1-s)(1-t)$ shows that $s=-t+t^2+ts+st-tst\in
 I$, whence $t\in I_q^\circ$ with quasi-adverse $s$.

If now $t\in\cl^\circ(I_q^\circ)$ and $a(1-t)+b=1$ for some $a,b$ in $R$, then
$xt-x-t+yb=0$ for some $x,y$ in $I$ by Lemma 4.6, whence $t-ryb\in
I_q^\circ$ for some $r$ in $I$. From the first part of the proof it
follows that $1-t+ryb\in R_q^{-1}$, which shows that $1-t\in\cl(R_q^{-1})$.  

Conversely, if $1-t\in\cl(R_q^{-1})$, and if $st-t-s+r=0$ for some
$s,r$ in $I$, then evidently $(1-s)(1-t)+r=1$ in $R$. By Lemma 4.6
this implies that $1-0=(1-x)(1-t)+yr$ for some $x,y$ in $I$. Since
$1-t\in\cl(R_q^{-1})$ it follows that $1-t+cyr\in R_q^{-1}$ for some
$c$ in $R$, whence $t-cyr\in I_q^\circ$ by the first part of the proof.
As $cy\in I$ this proves that $t\in\cl^\circ(I_q^\circ)$, as
desired. \hfill$\square$
\enddemo

\medskip

\example{4.8. Remark} We are now in a position to show that the notion
of $QB-$ring is also symmetric in the non-unital case. Indeed, let
$I$ be a non-unital ring and let $R=\widetilde{I}$
(as in 4.1). We identify $I$ with a two-sided ideal of $R$. It is
enough to see that $\clr^\circ(I_q^\circ)=I$ provided that
$\cl^\circ(I_q^\circ)=I$. Let $t$ in $I$, and assume that $tx-t-x+b=0$ 
for some $x,b$ in $I$. In $R$ this reads as $(1-t)(1-x)+b=1$. By the
right-handed version of Lemma 4.6 there are elements $x_1$ and $y$ in 
$I$ such that $(1-t)(1-x_1)+by=1$. Since $x_1\in\cl^\circ(I^\circ_q)$ 
it follows by Proposition 4.7 that $1-x_1\in\cl(R^{-1}_q)$. Now, by
Lemma 3.5, there exists $z$ in $R$ such that $1-t+byz\in
R_q^{-1}$. Again Proposition 4.7 implies that $t-b(yz)\in
I_q^\circ$. Since $yz\in I$, we see that $t\in\clr^\circ(I_q^\circ)$.
\endexample

\medskip

\proclaim{4.9. Theorem} An ideal $I$ in a unital ring $R$ is a $QB-$ring if
and only if
$$
1-I\subset\cl(R_q^{-1})\,.
$$
\endproclaim

\demo{Proof} By Proposition 4.7 the condition $1-I\subset\cl(R_q^{-1})$ is
equivalent to $I\subset\cl^\circ(I_q^\circ)$, which by definition means that
$I$ is a $QB-$ring \hfill$\square$ 
\enddemo

\medskip

\proclaim{4.10. Corollary} Every ideal in a $QB-$ring (unital
or not) will be a $QB-$ring. 
\endproclaim

\demo{Proof} Only the non-unital case deserves an argument. But if $I$ is an
ideal in the non-unital $QB-$ring $R$, then $I$ is also an ideal of $\wtR$, so
by Proposition 4.7
$$
\cl^\circ(I^{\circ}_q)=(1-\cl(\wtR^{-1}_q))\cap I = \cl^\circ(R^{\circ}_q)\cap I=R\cap I
=I\;.
$$
\hfill$\square$\enddemo

\example{4.11. Example} The ring $\Bbb Z[[x]]$ of formal power series in one
variable over $\Bbb Z$ is not a $B-$ring, but the ideal $x\Bbb Z[[x]]$ is. To
see this, recall first that $\Bbb Z[[x]]$ is well supplied with units, in fact
$a=\sum \alpha_nx^n \in (\Bbb Z[[x]])^{-1}$ if and only if $\alpha_0=\pm
1$. It follows that for any equation
$$
(1-y)(1-a)+b=1\,,
$$
with $a, b$ and $y$ in $x\Bbb Z[[x]]$, we can find an element $z$ in $x\Bbb
Z[[x]]$, such that $1-a+zb \in (\Bbb Z[[x]])^{-1}$. Indeed, $z=0$ will do.

We claim that there is no unital $B-$ring or even $QB-$ring $R$ that contains
$x\Bbb Z[[x]]$ as an ideal, and we mention this fact to show that our specific
$QB-$defini\-tions are needed for the non-unital case.

To establish the claim, assume that $R$ was such a $QB-$ring. The trivial
equation $3\cdot 2-5=1$ in $R$ shows that $2-5y\in R^{-1}_q$ for some $y$ in
$R$, so that 
$$
(1-(2-5y)v)R(1-v(2-5y))=0
$$
for some $v$ in $R$. Multiplying left and right with $x$, and replacing $R$
with $x^2$ we obtain the equation 
$$
(x^2-(2x-5xy)vx)(x^2-xv(2x-5yx))=0
$$
in $x\Bbb Z[[x]]$. However, this is a prime ring, so one of the factors must
be zero, say the one to the left. Write $xy=\sum \alpha_nx^n$ and $vx=\sum
\beta_nx^n$, and note that the equation 
$$
x^2=(2x-5xy)vx = (2x - 5\sum_{n=1}^\infty \alpha_nx^n)(\sum_{n=1}^\infty
\beta_n x^n) 
$$
forces $1=(2-5\alpha_1)\beta_1$ in $\Bbb Z$, which is plainly
impossible. \hfill $\square$
\endexample

%\newpage
\vskip1truecm

\subhead\nofrills{\bigbf 5. Skew Corners in QB--Rings}\endsubhead

\bigskip

\definition{5.1. Definition} Let $p$ and $q$ be idempotents in a
ring $R$ such that $pRq\ne 0$. We say that an element $x$
in $pRq$ is {\it quasi-invertible}, and write $x\in(pRq)_q^{-1}$, if
there exist elements $a,b$ in $qRp$ such that
$$
(p-xa)\perp (q-bx)\;.
$$
As in 2.1 we can take $y=a+b-axb$ in $qRp$, to obtain the equations
$x=xyx$ and
$$
(p-xy)\perp (q-yx)\;.
$$
In particular, $qRp \ne 0$ and $y\in(qRp)_q^{-1}$. Replacing if necessary $y$ by
$yxy$ we may also assume that $y=yxy$, so that $x$ and $y$ are von
Neumann regular elements and partial inverses for each other.

\enddefinition

\medskip

\definition{5.2. Definitions} With $p,q$ and $R$ as in 5.1 we define
the subset $\cl\sptilde((pRq)_q^{-1})$ to be the set of elements $a$ in $pRq$
such that whenever $xa+b=q$ for some elements $x$ in $qRp$ and $b$ in
$qRq$ we have $a+yb\in(pRq)_q^{-1}$ for some $y$ in $pRq$. Symmetrically 
we define the subset $\text{\rm cr}\sptilde((pRq)_q^{-1})$ as the elements $a$ in
$pRq$ such that whenever $ax+b=p$ for some $x$ in $qRp$ and $b$ in $pRp$ then
$a+by\in(pRq)_q^{-1}$ for some $y$ in $pRq$.

We leave it to the reader to verify, using the computations in Lemma
3.5, that if $ax+b=q$ for some $a$ in $qRp, \,b$ in $qRq$ and $x$ in
$\cl\sptilde((pRq)^{-1}_q)$, then $a+by\in (qRp)_q^{-1}$ for some $y$ in
$qRp$. In particular,  $\cl\sptilde((pRq)_q^{-1})=pRq$ if and only if
$\text{\rm cr}\sptilde((qRp)_q^{-1})=qRp$  

We shall be exclusively concerned with the case where both
$\cl\sptilde((pRq)_q^{-1})=pRq$ {\it and} $\text{\rm cr}\sptilde((pRq)_q^{-1})=pRq$,
and will refer to this situation as $pRq$ being a {\it $QB-$corner}. It
follows that in this case we also have that $qRp$ is a $QB-$corner. 

For an idempotent $p$ in $R$ we now seemingly have two notions of
$QB-$structure, one regarding $pRp$ as a unital ring in its own right, and one
regarding it as a corner in $R$. Fortunately these coincide, as will be seen
from Theorem 5.5, cf\. Corollary 5.7.
\enddefinition

\smallskip

As already mentioned in Theorem 2.3 we say that two idempotents are 
{\it(Murray-von Neumann) equivalent}, in symbols $p\sim q$, if $p=uv$ and
$q=vu$ for some elements $u,v$ in $R$. Replacing if necessary $u$ and $v$ with
$uvu$ and $vuv$, we may assume that $u$ and $v$ are von Neumann regular
elements and partial inverses for each other, and that
$$
(1-p)u=u(1-q)=0=(1-q)v=v(1-p)\;.
$$
We shall tacitly make these assumptions when we write $p\sim q$.

\medskip

\proclaim{5.3. Lemma {\rm (Cf\.\, \cite{{\bf 8}, 2.1})}} Let $p$ and $q$
be idempotents in a unital ring $R$ such that $1-p\sim 1-q$, i.e.
$1-p=uv$ and $1-q=vu$. If $x\in Rq$ such that $u+x$ is left or right
invertible in $R$, then the same holds for any element of the form
$u+x+(1-p)yq$. In particular for $u+px$.
 \endproclaim

\demo{Proof} Using the matrix decompositions
$R=(p+(1-p))R\left(\smallmatrix q\\ 1-q\endsmallmatrix\right)$ for $u+x$ and
$R=(q+(1-q))R\left(\smallmatrix p\\ 1-p\endsmallmatrix\right)$ 
for its inverses we can write $u+x=\left(\smallmatrix x_1 & 0\\ 
x_2 & u\endsmallmatrix\right)$.
If this element has a left inverse $a$ we have a matrix equation 
 $$
\pmatrix
a_{11} & a_{12}\\
a_{21} & a_{22}\endpmatrix\;
\pmatrix
x_1 & 0\\
x_2 & u\endpmatrix=
\pmatrix
q & 0\\
0 & 1-q\endpmatrix\;.
 $$
 Thus, $a_{22}u=1-q$, whence
 $$
a_{22}=a_{22}(1-p)=a_{22}uv= (1-q)v=v\;.
 $$
Moreover, $a_{12}u=0$, whence 
$$
a_{12}=a_{12}(1-p)=a_{12}uv=0\;.
$$
 This means that
 $$
q=a_{11}x_1+a_{12}x_2=a_{11}x_1\;.
 $$
 Thus, for any $y$ in $(1-p)Rq$ we have
 $$
\pmatrix
a_{11} & 0\\
-vya_{11} & v\endpmatrix\;
\pmatrix x_1 & 0\\
y & u\endpmatrix=
\pmatrix
q & 0\\
0 & 1-q\endpmatrix\;,
 $$
 so $u+px+y$ is left invertible.

Assume now that $u+x$ has a right inverse, i.e.
 $$
\pmatrix
x_1 & 0\\
x_2 & u\endpmatrix\;
\pmatrix
a_{11} & a_{12}\\
a_{21} & a_{22}\endpmatrix =
\pmatrix
p & 0\\
0 & 1-p\endpmatrix\;.
 $$
 This implies that $x_1a_{11}=p$, so for any $y$ in $(1-p)Rq$ we have
 $$
\pmatrix
x_1 & 0\\
y & u\endpmatrix\;
\pmatrix
a_{11} & 0\\
-vya_{11} & v\endpmatrix=
\pmatrix
p & 0\\
0 & 1-p
\endpmatrix\;,
 $$
so that $u+px+y$ is right invertible. \hfill$\square$
\enddemo

\medskip

\proclaim{5.4. Corollary} If $p$ and $q$ are idempotents in a unital
ring $R$ with $1-p\sim 1-q$, i.e. $1-p=uv$ and $1-q=vu$, and if
$u+x\in R_q^{-1}$ for some $x$ in $Rq$, then for any $y$ in $(1-p)Rq$
we also have $u+x+y\in R_q^{-1}$. In particular, $u+px\in R_q^{-1}$. 
\hfill $\square$
\endproclaim

\medskip

\proclaim{5.5. Theorem} Let $p$ and $q$ be idempotents in a unital
ring $R$ such that $1-p=uv$ and $1-q=vu$, and consider a non-zero
element $x$ in $pRq$.
Then $x\in (pRq)^{-1}_q$ if and only if $u+x \in R^{-1}_q$. Moreover, 
$x\in\cl\sptilde ((pRq)^{-1}_q)$ if and only if $u+x \in \cl (R^{-1}_q)$. 
\endproclaim

\demo{Proof} If $x\in (pRq)^{-1}_q$ we can find an element $y$ in $qRp$ such
that $p-xy \perp q-yx$. Set $u'= u+x$ and $v'=v+y$. Then by computation
$1-u'v'=1-uv-xy=p-xy$, and similarly $1-v'u'=q-yx$, so that $1-u'v'\perp
1-v'u'$ and $u+x\in R^{-1}_q$.

Conversely, if $u'=u+x\in R^{-1}_q$ and $v'$ is a quasi-inverse for $u'$ in
$R$, so that $1-u'v'\perp 1-v'u'$, then evidently $p-pu'v'p\perp
q-qv'u'q$. Since $pu'=x=u'q$ we can use $y=qv'p$ in $qRp$ to write $p-xy \perp
q-yx$, i\.e\. $x\in (pRq)^{-1}_q$.

Assume now that that $x\in \cl\sptilde ((pRq)^{-1}_q)$, and with $a=u+x$
consider an equation $za+b=1$ in $R$. Multiplying left and right with $q$ this
leads to the equation $qzpx+qbq=q$, since $aq=x=px$. By assumption we
therefore have $x+ybq\in (pRq)^{-1}_q$ for some $y$ in $pRq$. From what we
proved above this implies that $u+x+ybq\in R^{-1}_q$. However, $x+ybq \in pR$, so by
(the symmetric version of) Corollary 5.4 it follows that also $u+x+yb \in
R^{-1}_q$. Consequently $u+x\in \cl(R^{-1}_q)$.

Finally assume that $u+x\in \cl (R^{-1}_q)$, and consider an equation $zx+b=q$
for some elements $z$ in $qRp$ and $b$ in $qRq$. Then
$$
(v+z)(u+x)+b=1-q+zx+b=1\;.
$$
Since $u+x\in\cl(R_q^{-1})$ this implies that
$$
u+x+yb\in R_q^{-1}\;,\quad y\in R\;.
$$
As $x+yb\in Rq$ we conclude from Corollary 5.4 that also $u+x+pyb\in
R_q^{-1}$, i.e. we may assume that $y\in pRq$. By definition there is
an element $c$ in $R$ with
$$
(1-(u+x+yb)c)\perp (1-c(u+x+yb))\;.
$$
Multiplying the two elements  from left and right with $p$ and with $q$,
respectively, we see that
$$
(p-(x+yb)qcp) \perp (q-qcp(x+yb))\;,
$$
which shows that $x+yb\in(pRq)_q^{-1}$, whence
$x\in\cl\sptilde((pRq)_q^{-1})$. \hfill$\square$ 
\enddemo

\medskip

\proclaim{5.6. Corollary} If $p$ and $q$ are idempotents in a unital ring $R$
such that $1-p=uv$ and $1-q=vu$, and if $pRq\ne 0$, then $pRq$ is a
$QB-$corner if and only if\newline
\newline

\qquad\qquad\qquad\qquad $u+pRq\;\subset\;\cl(R_q^{-1})\;\bigcap\;\text{\rm cr}(R_q^{-1})\,.
\hfill\qed$ 
\endproclaim

\medskip

\proclaim{5.7. Corollary} For any non-zero idempotent $p$ in $R$ the ring
$pRp$ is a $QB-$corner, hence also a $QB-$ring, if and only if $1-p +
pRp\subset\cl(R_q^{-1})$.  \hfill$\square$ 
\endproclaim

\medskip

\proclaim{5.8. Corollary} If $R$ is a $QB-$ring, then $pRp$ is
a $QB-$corner for any non-zero idempotent $p$ in $R$, and  for any pair of
idempotents $p$ and $q$ such that $1-p\sim 1-q$ we either have $p\perp q$ or
else $pRq$ is a $QB-$corner. \hfill $\square$ 
\endproclaim

\medskip
 
Using the preceding results we can take up again the discussion from \S 2
about the r\^ole of $R^{-1}_q$ as the maximally extended von Neumann regular
elements. In $QB-$rings the solution is optimal:

\medskip
\proclaim{5.9. Theorem \text{\rm (Cf\. \cite{{\bf 8}, 2.6})}} Let $a$ be a
von Neumann regular element in a unital ring $R$. If
$a\in\cl(R_q^{-1})$, then $a\le u$ for some $u$ in $R_q^{-1}$.
\endproclaim

\demo{ Proof} By assumption $a=axa$ and $x=xax$ for some $x$ in
$R$. Define the idempotents $p=1-ax$ and $q=1-xa$, where $pa=aq=0$.
Thus, $xa+q=1$, and since $a\in\cl(R_q^{-1})$ this means that 
$a+yq\in R_q^{-1}$ for some $y$ in $R$. However, $1-p\sim 1-q$, so we can apply
Corollary 5.4 to conclude that also $a+pyq\in R_q^{-1}$.

If we now define $u=a+pyq$ then it follows from Lemma 2.7 that $a\le u$, as
desired.\hfill$\square$ 
\enddemo

\medskip

\proclaim{5.10. Corollary \text{\rm (Cf\. \cite{{\bf 8}, 4.3})}} In a unital 
$QB-$ring every von Neumann regular element extends to a quasi-invertible
element. \hfill $\square$
\endproclaim

\medskip

\proclaim{5.11. Corollary} In a unital $QB-$ring $R$ an element  
in $R^r$ is maximal with respect to the order $\le$ if and only if it belongs
to $R_q^{-1}$. 
\endproclaim

\demo{Proof} In view of Proposition 2.8 we need only consider an
element $a$ in $R^r$ which is maximally extended, and prove that
$a\in R_q^{-1}$. But that is evident from Corollary 5.10. \hfill$\square$ 
\enddemo

\medskip

\subhead{5.12. Morita Contexts}\endsubhead The concrete concept of a skew
corner $pRq$ in a ring $R$ and its transposed corner $pRq$ can be developed
abstractly as a theory of bimodules in Morita contexts.

If $R$ and $S$ are unital rings and if $M$ is an $R-S-$bimodule and
$N$ is an $S-R-$bimodule, we say that the pair $M,N$ is in a {\it
  Morita context} if we can find surjective bimodule maps
 $$
\varphi\colon M\otimes_S N\to R \quad\text{and}\quad \psi\colon N\otimes_R M\to S\;.
 $$
%i\.e\.
%$$
%a\varphi(x\otimes y)b=\varphi (ax \otimes yb)\quad\text{and}\quad c\psi(y
%\otimes x)d=\psi (cy \otimes xd) 
% $$
% for all $a,b$ in $R$, $c,d$ in $S$ and $x$ in $M$, $y$ in $N$.
Moreover these maps should satisfy the compatibility relations
$$
\varphi(x\otimes y)x'= x\psi(y\otimes x')\quad\text{and}\quad
y'\varphi(x\otimes y)=\psi(y'\otimes x)y
$$
for all $x, x'$ in $M$ and $y, y'$ in $N$.

Given a pair $M, N$ of bimodules in a Morita context for the rings $R$ and
$S$ we can define the {\it Morita ring} (or {\it linking ring})
$L=\left(\smallmatrix R & M\\ N & S\endsmallmatrix\right)$  with
pointwise addition and \lq\lq 
matrix product\rq\rq\; given by
 $$
\pmatrix
a_1 & x_1\\
y_1 & b_1\endpmatrix\;
\pmatrix a_2 & x_2\\
y_2 & b_2\endpmatrix=
\pmatrix
a_1a_2+\varphi(x_1 \otimes y_2) & a_1x_2+x_1b_2\\
y_1a_2+b_1y_2 & \psi(y_1 \otimes x_2)+b_1b_2
\endpmatrix\;.
 $$
 In $L$ we find the orthogonal pair of idempotents
 $$
p=\pmatrix
1_R & 0\\
0 & 0\endpmatrix\quad\text{and}\quad
q=\pmatrix
0 & 0\\
0 & 1_S\endpmatrix
 $$
 and the identifications
 $$
R=pLp\,,\;S=qLq\,,\;M=pLq\,,\;N=qLp\;.
 $$

The conditions that $\varphi$ and $\psi$ are surjective prove that $R$ and $S$
are Morita equivalent. Without this restriction we see that
$\varphi(M\otimes_S N)$ and $\psi( N \otimes_R M)$ will generate ideals $R_0$
and $S_0$ in $R$ and $S$, respectively, such that $M$ and $N$ will be in a
Morita context for $R_0$ and $S_0$. We can still form the large Morita
ring $L$, but we can also form the smaller Morita ring $L_0$, using
$R_0$ and $S_0$, which will then be an ideal in $L$ with quotient
$L/L_0=R/R_0\oplus S/S_0\,.$
 
\vskip1truecm
%\vskip0.5truecm
%\newpage

\subhead\nofrills{\bigbf 6. Matrices over QB--Rings}\endsubhead

\bigskip

We shall prove, with considerable effort, that $\Bbb M_n(R)$ is a $QB-$ring
whenever $R$ is. Our argument is an amalgamation of an argument attributed to
Kaplansky for the fact that Bass stable rank one passes to matrices, and the
series of reductions found in the proof of \cite{{\bf 8}, Theorem 4.5}.

First we show that \cite{{\bf 8}, Proposition 4.4} generalizes to the purely
algebraic setting.

\proclaim{6.1. Lemma} Consider elements $u,v$ in $R_q^{-1}$ with
quasi-inverses $x,y$ so that we have $(1-ux)\perp (1-xu)$ and $(1-vy)\perp
(1-yv)$. If moreover   
$$
uv+(1-ux)R(1-yv)\subset\cl(R_q^{-1})\cap\text{\rm cr}(R_q^{-1})\,,\tag $*$
$$ 
in particular if $R$ is a $QB-$ring, then either $(1-ux)R(1-yv)=0$ or it is a
$QB-$corner.  
\endproclaim

\demo{Proof} Let $I$ denote the ideal generated by the two idempotents $1-xu$
and $1-vy$. Since $((1-ux)R(1-yv))\cap I=0$ it follows that $(1-ux)R(1-yv)$
is isomorphic to its image in $R/I$. Moreover, by Proposition 3.7 condition
$(*)$ is still valid in $R/I$. To establish the lemma we may
therefore assume that $I=0$, i.e. $1=xu=vy$. Consequently
$$
(uv)(yx)=ux\quad\text{and}\quad (yx)(uv)=yv\;,
$$
so $ux\sim yv$. We can therefore apply Corollary 5.6 to conclude that
$(1-ux)R(1-yv)$ is a $QB-$corner if it is non-zero. \hfill$\square$
\enddemo

\medskip

\example{6.2. Remark} It is perhaps instructive to realize that
Lemma 6.1 also has a non-unital version: Consider elements $a,b$ in
$R_q^\circ$ with quasi-adverses $x,y$ so that
 $$
(ax-a-x)\perp (xa-a-x)\quad\text{and}\quad(by-b-y)\perp (yb-b-y)\,.
$$
If we can show that 
$$
ab-a-b + (ax-a-x)R(yb-b-y) \subset\cl^\circ(R_q^\circ)\cap\text{\rm cr}^\circ(R_q^\circ)\;,
$$
in particular if $R$ is a $QB-$ring, then $(ax-a-x)R(yb-b-y)$ is a $QB-$corner
if it is non-zero.

The proof is a tedious check of the fact that when we replace $a,b,x$ and $y$ by
$1-a$, $1-b$, $1-x$ and $1-y$ in $\wtR$, we can still apply Corollary 5.6.
\endexample

\medskip

To facilitate the arguments in the following we say that a right
unimodular row $(a,b)$ in $R$ (i.e. $aR+bR =R$) is {\it right reducible} if
$a+by\in R_q^{-1}$ for some $y$ in $R$. Note that we deliberately
choose a non-symmetric version, favouring the first coordinate in the
row. In the applications we shall consider a pair $(a,b)$ such that
$ax+b\in R^{-1}$, and right reducibility therefore means that $a\in\text{\rm
cr}(R_q^{-1})$.

The next lemma is a special case of the fact that multiplication by invertible
matrices preserve unimodular rows.

\proclaim{6.3. Lemma} Let $(a,b)$ be a right unimodular row in a
unital ring $R$. Then for any pair of units $u,v$ in
$R^{-1}$ and $c$ in $R$ the pair $(vau+vbc,vb)$ is also right
unimodular, and it is right reducible if and only if $(a,b)$ is right reducible.
 \endproclaim

\medskip

\demo{Proof} If $ax+by=1$ then
$(vau+vbc)u^{-1}xv^{-1}+vb(y-cu^{-1}x)v^{-1}=1$, so the new row is
still right unimodular.

Assume now that $(a,b)$ is right reducible, more specifically, $a+by\in
R_q^{-1}$. Then $vau+vbc+vb(yu-c)=v(a+by)u\in R_q^{-1}$, as
desired. The converse implications are similar.\hfill$\square$
\enddemo

\medskip

 \proclaim{6.4. Theorem} If $R$ is a unital $QB-$ring, then 
$\Bbb M_n(R)$ is also a $QB-$ring for any natural number $n$.
 \endproclaim

\demo{Proof} It suffices to prove the theorem for $n=2$, since by
iteration this would give the result for all numbers $n=2^k$; and
since $\Bbb M_n(R)$ is a corner in $\Bbb M_{2^k}(R)$ for $n\le 2^k$,
the general case follows from Corollary 5.8.

Consider therefore a right unimodular row $(\left(\smallmatrix a & b\\
c & d\endsmallmatrix\right)$, $\left(\smallmatrix e & f\\ g &
h\endsmallmatrix\right))$ and assume that
 $$
\pmatrix
a & b\\
c & d
\endpmatrix\;
\pmatrix
a' & b'\\
c' & d'
\endpmatrix+\pmatrix
e & f\\
g & h
\endpmatrix=
\pmatrix
1 & 0\\
0 & 1
\endpmatrix\;.\tag{$*$}
 $$
 Since $aa'+bc'+e=1$ and $R$ is a $QB-$ring, we have $a+(bc'+e)z_1\in
R_q^{-1}$ for some $z_1$ in $R$. By Lemma 6.3 our original row is
right reducible if and only if this is so for the row with elements
 $$
\pmatrix
a & b\\
c & d
\endpmatrix\;
\pmatrix
1 & 0\\
c'z_1 & 1
\endpmatrix+
\pmatrix
e & f\\
g & h
\endpmatrix\;
\pmatrix
z_1 & 0\\
0 & 0\endpmatrix\qquad\text{and}\qquad
\pmatrix
e & f\\
g & h
\endpmatrix\;;
$$
and now the first matrix has the $(1,1)$-element $a+bc'z_1+ez_1$,
which is quasi-invertible, whereas its second column is unchanged.
Without loss of generality we may therefore consider the equation
$(*)$ under the assumption that $a\in R_q^{-1}$.

Choose a quasi-inverse $x$ for $a$, and apply Lemma 6.3 to obtain the 
new right unimodular row of matrices with elements
 $$
\pmatrix
1 & 0\\ -cx & 1\endpmatrix\;
\pmatrix
a & b\\ c & d\endpmatrix\;
\pmatrix
1 & -xb\\ 0 & 1\endpmatrix\qquad\text{and}\qquad
\pmatrix
1 & 0\\ -cx & 1\endpmatrix\;
\pmatrix
e & f\\ g & h\endpmatrix\;.
 $$
 Now the first matrix has the form
 $$
\pmatrix
a & (1-ax)b\\
c(1-xa) & d_0\endpmatrix\;,
 $$
 so we may assume, without loss of generality, that we have equation
$(*)$ under the further restriction that $axb=0$ and $cxa=0$.

Computing the $(2,2)$-element in $(*)$ we obtain the equation
$cb'+dd'+h=1$, and since $R$ is a $QB-$ring this implies that
$d+(cb'+h)z_2\in R_q^{-1}$ for some $z_2$ in $R$. Using Lemma 6.3 we
pass to the right unimodular row with elements
 $$
\pmatrix
a & b\\ c & d\endpmatrix\;
\pmatrix
1 & b'z_2\\ 0 & 1\endpmatrix+
\pmatrix
e & f\\ g & h\endpmatrix\;
\pmatrix
0 & 0\\ 0 & z_2\endpmatrix\qquad\text{and}\qquad
\pmatrix
e & f\\ g & h\endpmatrix\;;
 $$
 where now the first matrix has the form
 $$
\pmatrix 
a & b+(ab'+f)z_2\\
c & d+(cb'+h)z_2
\endpmatrix\;,
 $$
 so its $(2,2)$-element is quasi-invertible. For the $(1,2)$-position
we compute
 $$
ax(b+(ab'+f)z_2)=0+ax(ab'+bd'+f)z_2=0
 $$
 using the $(1,2)$-equation from $(*)$ and the fact that $axb=0$.
Without loss of generality we may therefore assume that we have
equation $(*)$ under the further restriction that $d\in R_q^{-1}$.

Choose now a quasi-inverse $y$ for $d$ and apply Lemma 6.3 to obtain a 
new right unimodular row, where the first matrix has the form
 $$
\pmatrix
1 & -by\\ 0 & 1\endpmatrix\;
\pmatrix a & b\\ c & d\endpmatrix\;
\pmatrix
1 & 0\\ -yc & 1\endpmatrix=
\pmatrix
a & b(1-yd)\\
(1-dy)c & d\endpmatrix\;. 
 $$
 Without loss of generality we may therefore assume that we have
equation $(*)$ with the elements $a$ and $d$ in $R_q^{-1}$ with quasi-inverses
$x$ and $y$, such that
 $$
axb=0=byd\quad\text{and}\quad cxa=0=dyc\;.
 $$

The $(1,1)$-equation in $(*)$ followed by multiplication left and
right by $1-ax$, yields
 $$
(1-ax)b(1-yd)c'(1-ax)+(1-ax)e(1-ax)=1-ax\;.
 $$
Either $(1-ax)R(1-yd)=0$, in which case also the element $b$ is zero, and we do
nothing. Otherwise $(1-ax)R(1-yd)$ is a $QB-$corner by Lemma 6.1, whence 
$$
b+(1-ax)e(1-ax)z_3\in((1-ax)R(1-yd))_q^{-1}
$$
for some $z_3$ in $(1-ax)R(1-yd)$. Applying Lemma 6.3 we pass to the new right
unimodular row with elements 
 $$
\pmatrix
a & b\\ c & d\endpmatrix\;
\pmatrix
1 & -xez_3\\ 0 & 1\endpmatrix+
\pmatrix 
e & f\\ g & h\endpmatrix\;
\pmatrix
0 & z_3\\ 0 & 0\endpmatrix\qquad\text{and}\qquad
\pmatrix
e & f\\ g & h\endpmatrix\,;
 $$
 in which the first matrix has the form
 $$
\pmatrix
a & b+(1-ax)ez_3\\
c & d+g(1-ax)z_3\endpmatrix=
\pmatrix
a & b_1\\ c & d_1\endpmatrix\;,
$$
where $b_1\in((1-ax)R(1-yd))_q^{-1}$. Since $z_3\in R(1-yd)$ it
follows from Theorem 2.3 that $d_1\in R_q^{-1}$ with quasi-inverse
$y$ such that
$$
(1-d_1y)\perp (1-yd_1)\,,\quad (1-d_1y)\perp (1-yd)\,,\quad
(1-dy)\perp (1-yd_1)\;.
$$
Moreover,
$$
d_1y=(d+g(1-ax)z_3)y=dy
$$
since $z_3y=0$, so $d_1=d_1yd_1=dy(d+g(1-ax)z_3)\in d+dyR$.
Replacing $b$ and $d$ with $b_1$ and $d_1$ in $(*)$ we therefore
still have the equations
$$
d_1yc=0\quad\text{and}\quad axb_1=0\;.
$$
Furthermore,
$$
b_1yd_1=(b+(1-ax)ez_3)yd_1=byd_1=bydyd_1=0\;.
$$
By Lemma 6.3 we may therefore assume, on top of the previous
conditions, that we have ($*$) either with $b=0$ or with $b$ in 
$((1-ax)R(1-yd))_q^{-1}$.

A symmetric argument, using the $(2,2)$-equation in $(*)$ will
transform our right unimodular row to one with elements of the form
$$
\pmatrix
a_1 & b\\ c_1 & d\endpmatrix\qquad\text{and}\qquad
\pmatrix
e & f\\ g & h\endpmatrix\;,
$$
where now either $c_1=0$ and $a_1=a$ (in case $(1-dy)R(1-xa)=0$); or
else $a_1\in R_q^{-1}$ with
$xa_1=xa$, such that $a_1$ has $x$ as its quasi-inverse and
$$
c_1\in((1-dy)R(1-xa))_q^{-1}=((1-dy)R(1-xa_1))_q^{-1}\;.
$$
Applying Lemma 6.3 for the last time we may therefore assume that we have
equation $(*)$, but such that $a$ and $b$ are both quasi-invertible
with quasi-inverses $x$ and $y$ so that
$$
(1-ax)\perp (1-xa)\quad\text{and}\quad (1-by)\perp(1-yb)\;,
$$
and such that also
$$
b\in((1-ax)R(1-yd))_q^{-1} \cup \{ 0\} \qquad\text{and}\qquad
c\in ((1-dy)R(1-xa))_q^{-1}\cup \{ 0 \}
$$
with quasi-inverses $s$ in $(1-yd)R(1-ax)$ and $t$ in $(1-xa)R(1-dy)$ so that
$$
(1-ax-bs)\perp (1-yd-sb)\quad\text{and}\quad (1-dy-ct)\perp (1-xa-tc)\;.
$$
Straightforward computations show that
$$
\aligned
\pmatrix
a & b\\ c & d\endpmatrix\;
\pmatrix x & t\\ s & y\endpmatrix & =
\pmatrix
ax+bs & 0\\ 0 & dy+ct\endpmatrix\;,\\
\pmatrix
x & t\\ s & y\endpmatrix\;
\pmatrix
a & b\\ c & d\endpmatrix & =
\pmatrix
xa+tc & 0\\ 0 & yd+sb\endpmatrix\;,
\endaligned
$$
from which we deduce that $\left(\smallmatrix a & b\\ c &
d\endsmallmatrix\right)\in\Bbb M_2(R)_q^{-1}$ with quasi-inverse
$\left(\smallmatrix x & t\\ s & y\endsmallmatrix\right)$. Evidently
this means that the row is right reducible, and thus the original row from
$(*)$ is also right reducible, which proves that $\Bbb M_2(R) = \text{\rm
  cr}(\Bbb M_2(R))^{-1}_q$, as
desired. \hfill$\square$
\enddemo

\medskip

\example{6.5. Remark} Theorem 6.4 remains true also in the case
where $R$ is non-unital. Now instead of equation $(*)$ we must by Definition
4.4 consider an equation
 $$
\pmatrix
a & b\\ c & d\endpmatrix\;
\pmatrix
a' & b'\\ c' & d'\endpmatrix-
\pmatrix
a & b\\ c & d\endpmatrix-
\pmatrix
a' & b'\\ c' & d'\endpmatrix+
\pmatrix
e & f\\ g & h\endpmatrix=0\;,\tag{$**$}
$$
where all the matrices belong to $\Bbb M_2(R)$. Working in the
unital ring $\Bbb M_2(\wtR)=\Bbb M_2(R)+\Bbb M_2(\Bbb Z)$, but using
only matrices of the form $1-x$, where $x\in\Bbb M_2(R)$, we can
rewrite $(**)$ in the form $(*)$. Now all the matrix elements in
$(*)$ belong to $R$, except $a,a',d$ and $d'$, which are of the form
$1-R$, where, of course, $1$ denotes the unit in $\wtR$. Each of the reduction
steps in the proof of Theorem 6.4 will respect this structure. At the point
where we invoke Lemma 6.1 to transform $b$ (and later $c$) to a
quasi-invertible element in a skew corner, it is well to recall that
Lemma 6.1 has a non-unital version, cf. Remark 6.2.
\endexample

\medskip

\proclaim{6.6. Corollary} Let $R$ and $S$ be unital rings. If $R$ is a 
$QB-$ring and $S$ is Morita equivalent to $R$, then $S$ is also a
$QB-$ring.
\endproclaim
\demo{Proof} If $R$ and $S$ are Morita equivalent, then there is a
positive integer $n$ and an idempotent $e$ in $M_n(R)$ such that
$S\simeq eM_n(R)e$ (see, e.g., \cite{{\bf 1}, Corollary 22.7}). Then
the result follows using Theorem 6.4 and Corollary 5.8. \hfill$\square$
\enddemo

\subhead{6.7. Prime Embeddings}\endsubhead We say that a subring $S$
of a ring $R$ is {\it primely embedded} if $p\perp q$ in $S$ implies $p\perp
q$ in $R$ for any pair of idempotents $p$ and $q$ in $S$.

We have already used the fact that any ideal $I$ in a ring 
$R$ is primely embedded. Indeed, if $p\perp q$ in $I$ for two idempotents,
then $pRq=p^2Rq^2\subset pIq=0$, and similarly $qRp=0$, so that $p\perp q$ in
$R$.   

It is clear that if $R$ is a unital ring and $S$ is a primely embedded subring
containing  $1$, then $S_q^{-1}\subset R_q^{-1}$. Similarly  we have
$S_q^\circ\subset R_q^\circ$ in the non-unital case. 

\medskip

\proclaim{6.8. Proposition} Let $(R_n)$ be a sequence of $QB-$rings,  
and assume that we have homomorphisms $\varphi_n:R_n\to R_{n+1}$, such that 
$\varphi_n(R_n)$ is primely embedded in $R_{n+1}$ for every $n$. Then 
$R=\varinjlim R_n$ is a $QB-$ring.
\endproclaim

\demo{Proof} The elements in $R$ may be realized as (equivalence classes of)
sequences $x=(x_n)$, where $x_n\in R_n$ and $\varphi_n(x_n)=x_{n+1}$
eventually (i.e. from a certain $n_0$ onwards). (Two sequences being
equivalent if they agree eventually.) 

If $a=(a_n)\in R$ and $ax-a-x+b=0$ for some elements $x=(x_n)$ and
$b=(b_n)$ in $R$, then $\varphi_n(a_n)=a_{n+1}$,
$\varphi_n(b_n)=b_{n+1}$ and $\varphi_n(x_n)=x_{n+1}$ for all $n\ge
m$ for some $m$. Since $R_m$ is a $QB-$ring we can find $y_m$ in $R_m$ such
that $a_m-y_mb_m=c_m\in(R_m)_q^\circ$. Define $c=(c_n)$ for $n\ge m$
recursively by $c_{n+1} = \varphi_n(c_n)$. If $c_n\in(R_n)_q^\circ$ then
$\varphi_n(c_n)\in(\varphi_n(R_n))_q^\circ$ by Proposition 3.7. Since
$\varphi_n(R_n)$ is primely embedded in $R_{n+1}$ it follows that
$c_{n+1}\in(R_{n+1})_q^\circ$. As $c_m\in(R_m)_q^\circ$ we see by induction
that $c_n\in(R_n)_q^\circ$ for all $n$, whence $c\in R_q^\circ$. Defining
$y=(y_n)$ in $R$ starting with $y_m$ and inductively setting
$y_{n+1}=\varphi_n(y_n)$ we obtain the equation $a-yb=c$ in $R$,
whence $a\in\cl^\circ(R_q^\circ)$. Since $a$ was arbitrary, $R$ is a
$QB-$ring. \hfill$\square$ 
\enddemo

\medskip

\example{6.9. Remark} For each ring $R$ and each natural number $n$
there is a canonical embedding $\iota:\Bbb M_n(R)\to\Bbb M_{n+1}(R)$,
where $\iota(a)_{ij}=a_{ij}$ if $1\le i,j\le n$, but
$\iota(a)_{ij}=0$ if either $i=n+1$ or $j=n+1$. This embedding is not
unital, so the direct limit $\Bbb M_{\infty}(R)=\varinjlim\Bbb M_n(R)$
is a non-unital ring consisting of those matrices
$a=(a_{ij})$ over $R$ such that $a_{ij}=0$ if $i+j\ge m$ for some
$m$ (depending on $a$).
\endexample
\medskip

Combining Theorem 6.4 (maybe in its non-unital version described in
Remark 6.5) and Proposition 6.8 (which is allowed, since each
embedding $\iota$ certainly is prime) we therefore obtain the following
result: 

\medskip

\proclaim{6.10. Corollary} If $R$ is a $QB-$ring, then so is $\Bbb
M_{\infty}(R)$. \hfill $\square$ 
\endproclaim
%
%\medskip
%Recall that a ring $R$ is said to be {\it $\sigma-$unital} provided
%that there exists a sequence $\{a_n\}$ in $R$ such that
%$a_na_{n+1}=a_{n+1}a_n=a_n$ for all $n$, and for all $x$ in $R$, there 
%is $n_0$ in $\Bbb N$ such that $x=a_nx=xa_n$ for all $n\geq n_0$.
%\medskip
%\proclaim{6.10. Corollary} Let $R$ and $S$ be $\sigma-$unital rings. If 
%$R$ is a $QB-$ring and $S$ is Morita equivalent to $R$, then $S$ is
%also a $QB-$ring.
%\endproclaim
%\demo{Proof} If $R$ and $S$ are Morita equivalent, then the
%rings $\Bbb M_{\infty}(R)$ and $\Bbb M_{\infty}(S)$ are isomorphic, by
%\cite{{\bf 2}}. Then the result follows using Corollary 6.9. \hfill $\square$
%\enddemo

%\newpage
\vskip1truecm
%\vskip0.5truecm

\subhead\nofrills{\bigbf 7. Extensions of QB--Rings }
\endsubhead

\bigskip

\proclaim{7.1. Proposition} If $I$ is an ideal in a unital
$QB-$ring $R$, then
$$
(R_q^{-1}+I)/I=(R/I)_q^{-1}\;.
$$
\endproclaim

\demo{Proof} Let $\pi\colon R\to R/I$  denote the quotient morphism. From 
Proposition 3.7 we know that $\pi(R_q^{-1})\subset(R/I)_q^{-1}$, so it only
remains to show that every quasi-invertible element in $R/I$ lifts to a
quasi-invertible element in a $QB-$ring $R$.

For this we may consider $a$ and $b$ in $R$ such that $\pi(a)$ and $\pi(b)$
are elements in $(R/I)_q^{-1}$ with $(1-\pi(b)\pi(a))\perp (1-\pi(a)\pi(b))$
in $\pi(R)$. Thus, upstairs we have the relations
$$
(1-ba)R(1-ab)\subset I\quad\text{and}\quad (1-ab)R(1-ba)\subset I\;.
$$
Since $R$ is a $QB-$ring the trivial equation $ab+(1-ab)=1$ shows that we
have 
$$
v=b+y(1-ab)\in R_q^{-1}
$$
for some $y$ in $R$. Now choose a quasi-inverse $u$ for $v$ in $R^{-1}_q$. By
$(*)$ in Theorem 2.3 we then also have
$$
w=u+a(1-vu)+(1-uv)a\in R_q^{-1}\;.
$$
Passing to $R/I$ we evidently get $\pi(v)\pi(u)\pi(v)=\pi(v)$. But
we also have
$$
\aligned
& \pi(v)\pi(a)\pi(v)=\pi((b+y(1-ab))a(b+y(1-ab)))\\
& = \pi(ba(b+y(1-ab)))=\pi(b+y(1-ab))=\pi(v)\;,
\endaligned
$$
since $(1-\pi(ba))\pi(y)(1-\pi(ab))=0$. Thus, both $\pi(a)$ and
$\pi(u)$ are partial inverses for $\pi(v)$, whence by Theorem 2.3
$$
\pi(a)=\pi(u)+\pi(a)(1-\pi(vu))+(1-\pi(uv))\pi(a)\;.
$$
It follows that $\pi(a)=\pi(w)$, so that $w$ is the required lift of $\pi(a)$
in $R_q^{-1}$. \hfill$\square$ 
\enddemo

\medskip

\proclaim{7.2. Theorem \text{\rm (Cf\.\,\cite{{\bf 8}, 6.1})}} Let $I$ be an
ideal in a unital ring $R$. Then $R$ is a $QB-$ring if and only if
the following conditions are satisfied: 
\roster
\item"(i)" $R/I$ is a $QB-$ring ;
\item"(ii)" $(R_q^{-1}+I)/I=(R/I)_q^{-1}$;
\item"(iii)" $I+R_q^{-1}\subset\cl(R_q^{-1})$.
\endroster
\endproclaim

\demo{ Proof} If $R$ is a $QB-$ring then the first two conditions are
satisfied by Corollary 3.8 and Proposition 7.1, while
Condition (iii) is trivially true.

Conversely, if the three conditions are satisfied, take an arbitrary
element $a$ in $R$ and assume that $xa+b=1$ for some $x ,b$ in $R$. If
$\pi\colon R\to R/I$ denotes the quotient morphism, then 
$\pi(x)\pi(a)+\pi(b)=1$ in $R/I$, and since $R/I$ is a $QB-$ring this
implies that $\pi(a)+\pi(yb)\in(R/I)_q^{-1}$ for some $y$ in $R$. By
Condition (ii) we can find an element $t$ in $I$ such that $a+yb-t\in
R_q^{-1}$.

Using Condition (iii) this implies that $a+yb\in\cl(R_q^{-1})$, whence
$a\in\cl(\cl(R^{-1}_q))$ by the definition of $\cl$. But then
$a\in\cl(R_q^{-1})$ by Condition (iii) in Lemma 3.2, and thus
$R = \cl(R_q^{-1})$, as desired. \hfill$\square$
\enddemo

\medskip

\example{7.3. Remark} With obvious modifications the results in 7.1 and 7.2
remain true also in the non-unital case. If $R$ is a non-unital ring with an
ideal $I$ and quotient morphism $\pi:R\to R/I$, we consider the unitization
$\wtR=R\oplus\Bbb Z$ as in 3.1 and extend $\pi$ to a unital
morphism $\tilde\pi:\wtR\to R/I$ (if $R/I$ has a unit) or
$\tilde\pi:\wtR \to(R/I)^{\sim}$ (if $R/I$ has no unit) by setting
$\tilde\pi(1)=1$.

The argument in Proposition 7.1 now shows that $\pi(R_q^\circ)=(R/I)_q^\circ$
if $R$ is a $QB-$ring,  and in Theorem 7.2 we just have to replace the last
two conditions with
\roster
 \item"(ii')" $(R_q^\circ +I)/I=(R/I)_q^\circ$;
 \item"(iii')" $I+R_q^\circ\subset\cl^\circ(R_q^\circ)$.
 \endroster
\endexample

\medskip

\example{7.4. Remark} Condition (iii) in Theorem 7.2 is not easy
to verify directly, so it is fortunate that it is vacuously satisfied
in a number of interesting cases, cf\.\, Theorems 7.11 \& 7.14.

On the other hand, the condition is of independent interest and we shall
devote some attention to it. Note first that by Theorem 4.9 the condition
implies that the ideal is a $QB-$ring. Secondly observe from Condition (x)
in Lemma 3.2 that if an ideal $I$ satisfies Condition (iii), then actually
$I+\cl(R^{-1}_q)\subset \cl(R^{-1}_q)$. 

Having identified the extremally rich $C^*-$algebras with those $C^*-$algebras
which are $QB-$rings, cf\. Proposition 9.1, it follows from \cite{{\bf 8},
  Example 6.12} that we can not in Theorem 7.2 replace Condition (iii) with the
weaker condition: $I$ is a $QB-$ring. By necessity this means that the
extension theory for $QB-$rings is somewhat more complicated than that
governing $B-$rings and exchange rings. 
\endexample

\medskip

\proclaim{7.5. Lemma} If $I$ is a $QB-$ideal in a unital ring
$R$, then
$$
I+R^{-1}\subset\cl(R_q^{-1})\;.
$$
\endproclaim

\demo{ Proof} Take $a$ in $R^{-1}$ and $t$ in $I$ and assume that
$x(a-t)+b=1$ for some $x$ and $b$ in $R$. Then $xa(1-a^{-1}t)+b=1$, so
by Lemma 4.6 there are elements $r,r'$ in $I$ such that
$(1-r')(1-a^{-1}t)+rb=1$.

Since $I$ is a $QB-$ring $a^{-1}t\in \cl (I^{\circ}_q)$, whence 
$1-a^{-1}t\in \cl (R^{-1}_q)$ by Proposition 4.7, so that 
$a-t\in \cl (R^{-1}_q)$. \hfill$\square$
\enddemo

\medskip

\proclaim{7.6. Lemma} If $I$ is a $B-$ideal in a unital ring
$R$, then
$$
I+R_q^{-1}\subset\cl(R_q^{-1})\;.
$$
\endproclaim

\demo{ Proof} Take $u$ in $R_q^{-1}$ and $t$ in $I$, and assume
that
$$
x(u-t)+b=1\tag{$*$}
$$
for some $x$ and $b$ in $R$. Choose a quasi-inverse $v$ for $u$ so that the
two idempotents $p=uv$ and $q=vu$ satisfy $(1-p)\perp (1-q)$. Now rewrite the
equation $(*)$ as
$$
1=xu(1-vt)-x(1-p)t+b\;.
$$
Using Lemma 4.6 (with $vt$ and $b-x(1-p)t$
in place of $a$ and $b$) and that $I$ is a $B-$ideal, it follows that
$w_1=1-vt+s(b-x(1-p)t)\in R^{-1}$ for some $s$ in $I$ (and actually
$w_1-1\in I$). It follows that
$$
\aligned
uw_1 &=u-pt+usb-usx(1-p)t\\
     &=u-t+usb+(1-usx(1-p))(1-p)t\;.
\endaligned
$$
Since $(usx(1-p))^2=0$ the element $w_2=1+usx(1-p)$ is invertible with 
$w_2^{-1}=1-usx(1-p)$. Moreover, $w_2u=w_2^{-1}u=u$. Therefore, with
$t'=tw_1^{-1}$ we have 
$$
w_2(u-t+usb)w_1^{-1}=u-(1-p)t'\in R^{-1}_q
$$
by Theorem 2.3. It follows that also $u-t+usb\in R_q^{-1}$, so
$u-t\in\cl(R_q^{-1})$.\hfill$\square$  
\enddemo

\medskip

\example{7.7. Remark} Inspection of the preceding proof shows that if
$I$ is a $B-$ideal in $R$, then we also have the relations 
\roster
 \item"(i)" $I+R^{-1}\subset\cl(R^{-1})$;
 \item"(ii)" $I+R_r^{-1}\subset\cl(R_r^{-1})$;
 \item"(iii)" $I+R_\ell^{-1}\subset\cl(R_\ell^{-1})$.
 \endroster
For the first two it is actually easier to use a direct argument, but
the proof of relation (iii) needs the full force of the argument in
Lemma 7.6.
\endexample
 
\medskip

\proclaim{7.8. Proposition} In any unital ring $R$ there
is a largest ideal, denoted by $I_{qb}$, such that
$I_{qb}+R_q^{-1}\subset\cl(R_q^{-1})$. 
\endproclaim

\demo{ Proof} If $I_1$ and $I_2$ both satisfy the condition
$I_i+R_q^{-1}\subset\cl(R_q^{-1})$, then by Condition (x) in Lemma 3.2 they
also satisfy $I_i+\cl(R_q^{-1})\subset\cl(R_q^{-1})$ for $i=1,2$, whence
$$
I_1+I_2+R_q^{-1}\subset I_1+\cl(R_q^{-1})\subset\cl(R_q^{-1})\;.
$$
Therefore $I_{qb}$ is simply the sum of all the ideals that satisfy the
desired condition. \hfill$\square$
\enddemo

\proclaim{7.9. Proposition \text{\rm (Cf\.\,\cite{{\bf 11}, 2.14})}} If $R$ is a
unital ring such that $R$ is additively generated by its units
$(R=R^{-1}+R^{-1}+\cdots)$ then
$$
I_{qb}=\{x\in R\mid x+\cl(R_q^{-1})\subset\cl(R_q^{-1})\}\;.
$$
\endproclaim

\demo{ Proof} Obviously $I_{qb}$ is contained in the set $S$
defined by the right side of the equation, so by maximality it suffices to
show that $S$ an ideal. From the definition we see $S+S\subset S$ and also
$R^{-1}SR^{-1}=S$, because $R^{-1}\cl(R_q^{-1})R^{-1}=\cl(R_q^{-1})$.
Since every element $a$ in $R$ has a representation $a=\sum a_i$ with
$a_i$ in $R^{-1}$ it follows that $aS\subset S$, and $Sa\subset S$,
whence $S$ is an ideal. \hfill$\square$
\enddemo

\medskip

\example{7.10. Remark} Lemma 7.6 shows that the $QB-$ideal $I_{qb}$ defined
above contains every ideal in $R$ of Bass stable rank one. Since the sum of
$B-$ideals is again a $B-$ideal,  $I_{qb}$ therefore contains the maximal
$B-$ideal $I_{b}$ in $R$. Unfortunately the sum of $QB-$ideals may fail to be
a  $QB-$ideal, cf\.\,\cite{{\bf 8}, Example 6.12}, so we can not hope to
describe $I_{qb}$ as ``the maximal $QB-$ideal''. In $C^*-$algebra theory one
may instead characterize $I_{qb}$ as the ``largest well-behaved $QB-$ideal'',
in the sense that $I_{qb}$ is the maximal ideal such that $I_{qb} + B$ is a
$QB-$algebra for any $QB-$subalgebra $B$ of $A$ containing $1$, such that
$B^{-1}_q\subset A^{-1}_q$, cf\. \cite{{\bf 11}, Theorem 2.14}. The proof, however,
depends heavily on topological arguments.
\endexample

\medskip

\proclaim{7.11. Theorem \text{\rm (Cf\.\,\cite{{\bf 8}, 6.3})}} If $I$ is a
$B-$ideal in a unital ring $R$, then $R$ is a $QB-$ring if and
only if the following conditions are satisfied:
\roster
 \item"(i)" $R/I$ is a $QB-$ring;
 \item"(ii)" $(R_q^{-1} + I)/I=(R/I)_q^{-1}$.
 \endroster
 \endproclaim

\demo{ Proof} By Theorem 7.2 we only need to show that the two
conditions are sufficient. However, given the conditions it follows
from Lemma 7.6 that Condition (iii) in Theorem 7.2 is satisfied,
whence $R$ is a $QB-$ring. \hfill$\square$
\enddemo

\medskip

\proclaim{7.12. Corollary} If $I$ is a $B-$ideal in a unital ring $R$, and
$S$ is a $QB-$subring of $R$ containing $1$, such that $S$ is primely embedded
in $R$ and $R=I+S$, then $R$ is a $QB-$ring.
 \endproclaim

\demo{ Proof} Since $R/I$ is isomorphic to $S/I\cap S$ we know that $R/I$ is
a $QB-$ring by Corollary 3.8. By Theorem 7.11 we therefore only
need to verify that $(R_q^{-1}+ I)/I=(R/I)_q^{-1}$. However, if $v\in(R/I)_q^{-1}$
there is  by Proposition 7.1 an element $u$ in $S_q^{-1}$ such
that $u+I=v$. Since $S$ is primely embedded in $R$ this means that 
$u\in R_q^{-1}$ cf\. 6.7.\hfill$\square$
\enddemo

\medskip

\example{7.13. Remark} For each unital ring $R$ its Jacobson radical $\Cal J(R)$ is
a $B-$ideal, since $1-x\in R^{-1}$ for any $x$ in $\Cal J(R)$. Moreover,
any lift of a left, respectively right invertible element in $R/\Cal J(R)$
will be left, respectively right invertible in $R$. It follows from
Theorem 7.11 that if $R/\Cal J(R)$ is a prime ring, then $R$ is a $QB-$ring
if and only if $R/\Cal J(R)$ is a $QB-$ring.

\endexample

\medskip

\proclaim{7.14. Theorem \text{\rm (Cf\.\,\cite{{\bf 8}, 6.6})}} If $I$ is a
$QB-$ideal in a unital ring $R$, such that $R/I$ is a $B-$ring, 
then $R$ is a $QB-$ring provided that $(R^{-1}+I)/I=(R/I)^{-1}$.
 \endproclaim

\demo{ Proof} Take any $a$ in $R$ and assume that $xa+b=1$ for
some elements $x$ and $b$ in $R$. Then with $\pi\colon R\to R/I$ the
quotient morphism we also have $\pi(x)\pi(a)+\pi(b)=1$ in $R/I$, and
since $R/I$ is a $B-$ring this implies that
$\pi(a)+\pi(y)\pi(b)=\pi(u)\in(R/I)^{-1}$ for some $y$ and $u$ in
$R$. Since $\pi(R^{-1})=(R/I)^{-1}$ we may assume that $u\in R^{-1}$, so
that we have the equation
$
a-t+yb=u$ for some $t$ in $I$. By Lemma 7.5 this implies that
$a+yb\in\cl(R_q^{-1})$. Therefore
$a\in\cl(\cl(R_q^{-1}))=\cl(R_q^{-1})$ (see Condition (iii) in Lemma
3.2). \hfill$\square$
\enddemo

\medskip

\proclaim{7.15. Corollary} If $I$ is an $QB-$ideal in a unital
ring $R$, and $S$ is a $B-$subring of $R$ containing
$1$, such that $R=I+S$, then $R$ is a $QB-$ring.
\endproclaim

\demo{ Proof} Since $R/I$ is isomorphic to $S/I\cap S$ we know that $R/I$ is a
$B-$ring, because Bass stable rank one is preserved under quotients. To apply 
Theorem 7.14 we therefore need only to verify that invertibles lift from
$R/I$. But invertibles certainly lift from quotients of $B-$rings, so if 
$x\in(R/I)^{-1}$\; ($= (S/(I\cap S))^{-1})$ there is an element $u$ in
$S^{-1}$ such that $u+I=v$. As $S^{-1}\subset R^{-1}$ we conclude that
$\pi (R^{-1})=(R/I)^{-1}$, as desired. \hfill$\square$
\enddemo

\medskip

\example{7.16. Toeplitz-like Examples} We present an example of a unital
$QB-$ring $S$ which is von Neumann regular and is an extension of
two $B-$rings, but has Bass stable rank two. As we shall see, this
example can be thought of as an algebraic analogue of the Toeplitz
algebra, since it is generated, in a suitable sense, by a unilateral
shift. Other examples of non-regular (even not exchange) rings will be
given later (cf\. 8.9). Our construction is modelled after the example
given in \cite{\bf 15}. We provide some details for the convenience
of the reader.

\medskip

Let $\Bbb F$ be a countable field and let $t$ be an indeterminate. Let
$\Bbb F(t)$ be the field of rational functions on $\Bbb F$, and let
$\delta$ be the valuation on $\Bbb F(t)$ associated with the maximal ideal
$(t)$ of $\Bbb F[t]$, i\.e\. $\delta (0)=+\infty $ and $\delta
(t^nf(t)/g(t))=n$ where $t$ does not divide $f(t)g(t)$. Let $V=\{ x\in
\Bbb F(t)\mid \delta (x)\ge 0\}$ be the valuation ring associated with $\delta
$. Note that $V$ is a local ring with maximal ideal $\{ x\in \Bbb F(t)\mid
\delta (x)>0 \}$.

\medskip

We claim that the vector space $W=\Bbb F(t)$ has an $\Bbb F-$basis $\{ v_i
\}_{i\in \Bbb Z}$ 
such that $\delta (v_i)=i$ for all $i$ in $\Bbb Z$. First note that $W$ is
of countable dimension over $\Bbb F$ because $\Bbb F$ is
countable. Take an $\Bbb F-$basis $\{ w_n \} _{n\ge 
0}$ for $W$. We can modify this basis in order to get $\delta
(w_i)\ne \delta (w_j)$ for all $i\ne j$. In fact, if $\delta (w_n)=
\delta (w_i)$ for $i<n$ there is an element $\alpha _i$ in $\Bbb F$ such that
$w_n/w_i -\alpha _i\in tV$, which implies that $\delta (w_n-\alpha _iw_i)>
\delta (w_n)$. If $\delta (w_n-\alpha _iw_i)=\delta (w_j)$ for some
$j<n$, then the same argument shows that there is an element $\alpha
_j$ in $\Bbb F$
such that $\delta (w_n-\alpha _iw_i-\alpha _jw_j)>\delta (w_j)>\delta
(w_i)$. Thus we get an element $w_n-\alpha _iw_i-\cdots -\alpha_kw_k$,
such that $\delta (w_n-\alpha _iw_i-\cdots -\alpha _kw_k)\ne \delta
(w_t)$ for all $t<n$. Then we substitute $w_n$ by $w_n-\alpha
_iw_i-\cdots -\alpha_kw_k$.

Now assume that $\delta (w_i)\ne \delta (w_j)$ for $i\ne j$. Writing an
arbitrary element $v$ of $W$ as $v=\lambda_{i_1}w_{i_1}+\cdots +
\lambda _{i_k}w_{i_k}$ with $\delta (w_{i_1})<\delta (w_{i_2})<\cdots < \delta
(w_{i_k})$ , we have $\delta (v)=\delta (w_{i_1})$. Since $\delta
(t^i)=i$ for all $i$ in $\Bbb Z$ we have a 
bijective correspondence $\varphi :\Bbb Z\rightarrow \Bbb Z ^+$ such
that $\delta (w_{\varphi (i)})=i$. Take $v_i=w_{\varphi (i)}$ for all
$i$ in $\Bbb Z$. Note that $V=\langle \{v_i\}_{i\geq 0} \rangle$, and
obviously we may assume that $v_0=1$.

\medskip

Consider the representation $\lambda :\Bbb F(t)\rightarrow \End _\Bbb F (W)$
given by multiplication. Let $\pi:W\rightarrow W$ be the projection onto
$V$ with kernel $\langle \{v_i\}_{i<0}\rangle $. Identifying
$\pi\End _\Bbb F(W)\pi$ with $\End _\Bbb F (V)$, we may regard $\pi\lambda (\Bbb F(t))\pi$
as a subring of $\End _\Bbb F(V)$. Each endomorphism of $V$ has a
column-finite matrix associated to the basis $\{v_i\}_{i\geq 0}$ of $V$. Since
$\lambda (v)(v_i)=vv_i\in \langle v_{i+\delta (v) }, v_{i+\delta (v)+1}, \dots
\rangle $ (we use that $\delta(vv_i)=i+\delta (v)$), it follows that
$\pi \lambda (\Bbb F(t))\pi\subset  \Bbb B (\Bbb F)$, where $\Bbb B (\Bbb F)$ is
the algebra of row-and-column finite matrices over $\Bbb F$ (and we
identify an element of $\pi\lambda (\Bbb F(t))\pi $ with its matrix with respect
to the basis $\{ v_i\}_{i\geq 0}$ of $V$). Let $\Bbb M_\infty(\Bbb F)$ be 
the ideal of $\Bbb B (\Bbb F)$ consisting of matrices with only a finite
number of non-zero entries, cf\. 6.9, and consider the ring $S=\pi\lambda (\Bbb
F(t))\pi+\Bbb M_\infty(\Bbb F)$. Since $\pi\lambda (u)\pi\lambda (v)\pi-\pi\lambda
(uv)\pi\in \Bbb M_\infty(\Bbb F)$ for all $u,v$ in $\Bbb F(t)$, 
we see that $S$ is an $\Bbb F$-subalgebra of $\Bbb B (\Bbb F)$. There is a
surjective homomorphism $\rho :S\rightarrow \Bbb F(t)$ defined by $\rho
(\pi\lambda (u)\pi+m)=u$ for all $u$ in $\Bbb F(t)$ and all $m$ in $\Bbb
M_\infty(\Bbb F)$. Since $\Bbb M_\infty(\Bbb F)$ is 
regular and $\Bbb F(t)$ is also regular this shows that $S$ is a regular
ring.

Now let $a=\pi\lambda (t)\pi$ and $b=\pi\lambda (t^{-1})\pi$. Since $\lambda
(t)\pi=\pi\lambda (t)\pi$ we get $ba=1$. However, $ab$ is not equal to $1$,
and $1-ab$ is a one-dimensional idempotent in $\Bbb M_\infty(\Bbb F)$. It follows that
$S$ is not directly finite. In particular it is not a $B-$ring, and it 
is an extension of the $B-$rings $\Bbb M_\infty(\Bbb F)$ and $\Bbb F(t)$. Notice that $\bsr
(S)\leq \text{max}\{\bsr(\Bbb F(t)),\bsr(\Bbb M_\infty(\Bbb F))+1\}=2$, whence
$\bsr(S)=2$. Also we see that $a,b\in S_q^{-1}$, and since every non-zero
element in $\Bbb F(t)$ has the form $t^i v$, where $i\in \Bbb Z$ and $v$ is an
invertible in $V$, we get that every non-zero element in $\Bbb F(t)$ lifts to
a quasi-invertible element of $S$. (Note that the map $\pi\lambda (-)\pi $
provides an isomorphism from $V$ onto a subalgebra of $S$.)

It follows from Theorem 7.11 that $S$ is a $QB-$ring.
\endexample
\medskip

The next example is based on a construction due to Bergman (see
\cite{{\bf 16}, Example 5.10}, and also \cite{{\bf 20}}).

\example{7.17. Example} In the setting of Example 7.16, consider the
map $\rho:S\rightarrow \Bbb F(t)$. Let $S^o$ be the opposite ring of
$S$. Since $\Bbb F(t)$ is commutative we have an induced surjection
$\rho:S^o\rightarrow \Bbb F(t)$. Taking the pullback of both maps we get the
ring $T=\{ (x,\overline{y})\mid \rho (x)=\rho (y)\}$. The ring $T$ has a unique
maximal ideal, viz\. $\Bbb M_\infty(\Bbb F)\times \Bbb M_\infty(\Bbb F)^o$ and
$T/(\Bbb M_\infty(\Bbb F)\times \Bbb M_\infty(\Bbb F)^o) \simeq \Bbb F(t)$,
which proves that $T$ is a regular ring. The ring $T$ is directly finite by
\cite{{\bf 20}, Lemma 13}, but has a quotient isomorphic to $S$ and
a quotient isomorphic to $S^o$, so it is not a $B-$ring. The elements
$(a,\overline{a})$ and $(b,\overline{b})$ are quasi-invertible in $T$ (though they
are not right or left invertible) and it follows again from Theorem
7.11 that $T$ is a $QB-$ring.
\endexample

\vskip1truecm
%\vskip0.5truecm

\subhead\nofrills{\bigbf 8. Exchange Rings}\endsubhead 
\bigskip

\definition{8.1. Definitions} A unital ring $R$ is called an {\it exchange
ring} if for every element $a$ in $R$ there is an idempotent $p$ in $aR$ such
that $1-p\in(1-a)R$. This is not the original definition (which concerns a
finite exchange property in $R$-modules, cf. \cite{\bf 30}), but is an
equivalent description found by Goodearl and Nicholson, see \cite{\bf 17} and
\cite{\bf 21}. Rewriting the condition $1-p\in(1-a)R$ as an equation
$1-p=(1-a)(1-y)$ for some $y$ in $R$, i.e. $p=a+y-ay$, we obtain a definition
of an exchange ring suitable for the non-unital case, cf\. \cite{{\bf 2}}.

The class of exchange rings is extensive and includes all von Neumann
regular rings, all $\pi$-regular rings, the semi-perfect rings (which
are exactly the exchange rings that are semi-local), right
self-injective rings and $C^*-$algebras of real rank zero. (In fact, the
$C^*-$algebras  which are exchange rings are precisely those of real
rank zero, by \cite{{\bf 4}, Theorem 7.2}.)
The class of exchange rings is stable under ideals and quotients as well as
corners and matrix tensoring; and if $0\to I\to R\to Q\to 0$ is a
short exact sequence of rings, then $R$ is exchange if and only if
both $I$ and $Q$ are exchange and idempotents lift from $Q$ to $R$,
cf\. \cite {{\bf 2}, Theorem 2.2}. Evidently the class is also stable
under direct limits.
\enddefinition

\medskip

\proclaim{8.2. Proposition} In a unital, semi-primitive exchange
ring $R$ an element in $R^r$ is maximal with respect to the order
$\le$ if and only if it belongs to $R_q^{-1}$.
\endproclaim

\demo{Proof} By Proposition 2.8 we need only consider a maximally extended
element $a$ in $R^r$, and prove that  $a\in R_q^{-1}$.

Towards this end choose a partial inverse $x$ for $a$ and set $p=1-ax$ and
$q=1-xa$. If $pRq\ne 0$ we can find a non-zero element $y$ in $pRq$. Since $R$
is semi-primitive there is a maximal right ideal $J$ of $R$ such that $y\notin
J$. Hence $1-yc\notin R_r^{-1}$ for some $c$ in $R$. Since $R$ is an
exchange ring we can find an idempotent $e$ in $ycR$ such that
$1-e\in(1-yc)R$. We have made sure that $1-e\ne 1$, so $e\ne 0$. Since
$e=eycd$ for some $d$ in $R$, the element $ey\in pRq\setminus\{0\}$
and $(ey)(cd)(ey)=e^2ey=ey$ so $ey$ is von Neumann regular with $cd$
as partial inverse.

By Lemma 2.7 the element $b=a+ey$ is von Neumann regular in $R$ and properly
extends $a$, contradicting the maximality. Thus, after all,
$pRq=0$ (and also $qRp=0$), whence $a\in R_q^{-1}$. \hfill$\square$
\enddemo

\medskip

\example{8.3. Examples} Corollary 5.11 and Proposition 8.2 have different
implications, despite their similarity. In a $QB-$ring every von Neumann
regular element extends to a maximal one, but maybe they are all maximal to
begin with (except zero). This happens e\.g\. if $R=C([0,1])$, where each
non-zero regular element is invertible.

In an exchange ring, by contrast, there is an abundance of idempotents, and
therefore also a great variety of von Neumann regular elements. The problem
here is that they might not all extend to quasi-invertible elements. Of
course, for the idempotents there are no problems: each extends to 1. 

To construct an example of an exchange ring with a regular element that does
not extend, take $A=\Bbb B(\Cal H)\oplus \Bbb B(\Cal H)$. Represent $A$ with
infinite multiplicity on the Hilbert space $\Cal K=(\Cal H \oplus \Cal
H)\otimes \ell^2$ (so every operator in $A$ is repeated infinitely often along
the diagonal), and put $R= A + \Bbb K(\Cal K)$. Thus, $R$ is a split extension
of the algebra of compact operators on $\Cal K$ and $A$. Evidently both 
$\Bbb K(\Cal K)$ and $A$ are exchange rings (they are $C^*-$algebras of real rank
zero), and projections lift since the extension splits, so $R$ is an exchange
ring.  Let $s$ denote the unilateral shift in $\Bbb B(\Cal H)$ (or any other
non-unitary isometry). Then $u = s\oplus s^*\in A^{-1}_q$ with  quasi-inverse 
$u^* =s^* \oplus s$. However, since $R$ is a primitive algebra its quasi-invertible
elements are either left or right invertible, so $u\notin R^{-1}_q$. If $u$
could be extended to a left or right invertible element $w$ in $R$, then in the
quotient $R/\Bbb K(\Cal K)=A$ we would have $w-u+ \Bbb K(\Cal K) =0$, since
$u$ is maximally extended in $A$. But since both the kernel and the co-kernel
of $u$ on $\Cal K$ are infinite dimensional, no compact perturbation can make $u$
left or right invertible.

As the following result shows, the global obstruction to extension of regular
elements in an exchange ring $R$ is exactly that $R$ fails to be a $QB-$ring.
\endexample

\medskip

\proclaim{8.4. Theorem} If $R$ is a unital exchange ring the following
conditions are equivalent:
\roster
\item "(i)" $R$ is a $QB-$ring\,;
\item "(ii)" Every element in $R^r$ extends to an element in $R^{-1}_q$\,;
\item "(iii)" For every element $x$ in $R^r$ there is a $v$ in $R^{-1}_q$ such
  that $x=xvx$\,.
\endroster 

\endproclaim

\demo{Proof} (i) $\implies$ (ii) This is immediate from Corollary 5.10.

\noindent (ii)$\implies$(i) Given an equation $xa+b=1$ in $R$ we use the
exchange property to find an idempotent $p$ in $Rxa$, such that $1-p \in Rb$. 
Specifically, $p=rxa$ and $1-p = sb$. Then $ap$ is a von Neumann regular
element with partial inverse $rx$, and by assumption it extends to an element
$u$ in $R^{-1}_q$. Thus by definition, and using Lemma 2.6 if
necessary, $ap=u(rx)(ap)=up$. It follows that with
$y=(u-a)s$ we can write
$$
u=up+u(1-p)=ap+u(1-p)=a+(u-a)sb=a+yb\,,
$$
as desired.

\noindent(ii)$\iff$(iii) This is Proposition 2.9. \hfill $\square$
\enddemo

\medskip

The next lemma is an adaptation of \cite{{\bf 31}, Theorem 2.1} to our
situation.

\proclaim{8.5. Lemma} Let $R$ be a unital ring. Then the following
conditions are equivalent: \roster
\item "(i)" $R$ is a $QB-$ring\,;

\item "(ii)" Given any $R-$module $M$ and decompositions
  $M=A_1\oplus H=A_2\oplus K$ with 
$A_1\simeq A_2\simeq (R_R)^n$ for some $n\ge 1$, there exists a pair
of orthogonal ideals $I$ and $J$ and decompositions $M=E\oplus
B\oplus H=E\oplus C\oplus K$ such that $BI=B$ and $CJ=C$\,;

\item "(iii)" Given any $R-$module $M$ and decompositions
  $M=A_1\oplus H=A_2\oplus K$ with 
$A_1\simeq A_2\simeq (R_R)^n$ for some $n\ge 1$, there exists a pair
of orthogonal ideals $I$ and $J$ and decompositions $M=E_1\oplus
N=E_2\oplus N$ such that $A_1=E_1\oplus B $ and $A_2=E_2\oplus C$
with $BI=B$ and $CJ=C$\,.
\endroster
\endproclaim

\demo{Proof} (i)$\implies $(ii) For $i=1,2$ let $\rho'_i\colon M\to
A_i$ and $\rho_1\colon M\to H$ and $\rho _2\colon M\to K$ denote the
projections and $\tau'_i\colon A_i\to M$ and $\tau_1\colon H\to M$
and $\tau _2\colon K\to M$ the corresponding injections. Fix
isomorphisms $\iota_i\colon A_i\to R^n$ and put
$\pi_i=\iota_i\rho'_i$ and $\sigma_i=\tau'_i\iota_i^{-1}$. Then
$\pi_i\sigma_i=1_{R^n}$ and $\sigma_i\pi_i+\tau_i\rho_i=1_M$. It
follows that 
$$
1_{R^n}=\pi_2\sigma_2=\pi_2(\sigma_1\pi_1+\tau_1\rho_1)\sigma_2
=(\pi_2\sigma_1)(\pi_1\sigma_2)+(\pi_2\tau_1\rho_1\sigma_2)=ax+b\;.
$$ 
Identifying $\Bbb M_n(R)$ with $\End_R(R^n)$, and using that $\Bbb M_n(R)$
is a $QB-$ring, cf. Theorem 6.4, we can find $y$ in $\Bbb M_n(R)$ and $u$
in $(\Bbb M_n(R))^{-1}_q$, such that $u=a+by$, i\.e\. 
$$
u=\pi_2(\sigma_1+\tau_1\rho_1\sigma_2 y)\,. 
$$

Choose a quasi-inverse $v$ for $u$ and define the idempotents $p=vu$
and $q=uv$ in $\Bbb M_n(R)$. Moreover, let $\varphi = \sigma_1
+\tau_1\rho_1\sigma_2y$ in $\Hom(R^n,M)$. Then from the equation
above we see that $\pi_2\varphi=u$. We also compute
$\pi_1\varphi=\pi_1(\sigma_1 +\tau_1\rho_1\sigma_2y)=1_{R^n}$, since
$\pi_1\tau_1=0$.

Let $D_1=\ker p\pi_1$ and $D_2=\ker q\pi_2$, and put $E=\varphi
p(R^n)$ in $M$. If $m\in E\cap D_1$, then $m=\varphi p(a)$ for some
$a$ in $R^n$, where $$ 0=p\pi_1(m)=p\pi_1\varphi p(a)=
p1_{R^n}p(a)=p(a)\;. $$ Consequently $m=0$, and so $E\cap D_1=0$.
Similarly, $E\cap D_2=0$. Now take $m$ in $M$ and write $$ m= \varphi
p \pi_1 (m)+(m-\varphi p \pi_1(m))\;. $$ Then $p\pi_1(m-\varphi p
\pi_1 (m))= p\pi_1(m) - p 1_{R^n}p \pi_1(m)=0$, so $M=E\oplus D_1$.
Similarly $M=E\oplus D_2$.

Now we decompose further: Let $B=\sigma_1(1-p)(R^n)$ and
$C=\sigma_2(1-q)(R^n)$ in $M$. Then 
$$ 
D_1=\ker p\pi_1=\ker\pi_1\oplus \sigma_1 (1-p)(R^n) = H\oplus B\;. 
$$ 
Similarly, 
$$
D_2=\ker q\pi_2=\ker \pi_2\oplus \sigma_2 (1-q)(R^n) = K\oplus C\;.
$$ 
Finally, let $\Bbb M_n(I)$ and $\Bbb M_n(J)$ be the ideals in $\Bbb M_n(R)$
generated by $1-p$ and $1-q$, respectively, where $I$ and $J$ are
ideals in $R$. Since $(1-p)\perp (1-q)$ it follows that $I\perp J$,
and evidently $BI=B$ and $CJ=C$. Moreover, $E\oplus B\oplus H=E\oplus
D_1 =M=E\oplus D_2=E\oplus C\oplus K$, as desired.

(ii)$\implies $(i) Take $(a_1,a_2)$ in $R^2$ such that
$a_1R+a_2R=R$. Then $(a_1,a_2)$ provides a split epimorphism
$\varphi :M\rightarrow R$, where $M=R^2$, so we have decompositions
$M=R\oplus R= A_2\oplus K$, where $A_2\simeq R_R$ and $K$ is the
kernel of $\varphi$. Therefore there exists a pair of orthogonal
ideals $I$ and $J$ and a decomposition $$M=E\oplus B\oplus R=E\oplus
C\oplus K$$ with $BI=B$ and $CJ=C$. Write $E_1=\varphi (E)$ and
$C_1=\varphi (C)$. Then $R=E_1\oplus C_1$. Let $e$ in $R=\End(R_R)$
be the projection onto $E_1$ with kernel $C_1$, so that
$E_1=eR$. Notice that $1-e\in J$. Let $\psi _1 :R\rightarrow M$ be the
inverse of $\varphi$ restricted to $E\oplus C$ and let $\psi = \psi_1\circ
e$. Then $\psi $ is given by a pair $(x_1,x_2)$ in $M$ such that
$x_ie=x_i$ for $i=1,2$ and $a_1x_1+a_2x_2=e$. Let $\pi: M=R\oplus
R\rightarrow R$ be the projection onto the first factor. Write
$E_2=\pi (E)$ and $B_1=\pi (B)$. As before, $R=E_2\oplus B_1$. Let $f$
in $R=\End(R_R)$ be the projection onto $E_2$ with kernel $B_1$, so
that $E_2=fR$. Notice that $1-f\in I$. Now we have that $\pi\circ
\psi$ provides an isomorphism from $E_1$ onto $E_2$. Let $y_1$ in
$R=\End(R_R)$ be an element such that $y_1f=y_1=ey_1$ and implements
the inverse of the isomorphism $(\pi\circ \psi)\,|E_1$. Observe that
$y_1x_1=e$ and $x_1y_1=f$ so that $x_1$ and $y_1$ are quasi-invertible
and quasi-inverses for each other (because $1-e\in J$ and $1-f\in I$ and $I$  
and $J$ are orthogonal). Since $(a_1+a_2x_2y_1)x_1=y_1x_1=e$
this implies that $a_1+a_2x_2y_1$ is quasi-invertible (see, for
example, Theorem 2.3), showing that $R$ is a $QB-$ring.

The proof of (i)$\iff$(iii) is dual to the one of (i)$\iff$(ii), cf\.
\cite{{\bf 31}, Theorem 1.6}. Since we know that the notion of
$QB$-ring is symmetric (Theorem 3.6), we are done. \hfill$\square$
\enddemo
 
\medskip
 
\definition{8.6. Definitions} For a unital ring $R$ we denote by $\Cal V(R)$
the monoid of isomorphism classes of finitely generated, projective, right
$R-$modules. (The standard notation for the category of such modules
is $\Cal F\Cal P(R)$.) The addition in $\Cal V(R)$ is direct sum of
representatives, and $\Cal V(R)$ is ordered by $x\le y$ if $M\oplus P\simeq N$
for some representatives $[M]=x$ and $[N]=y$. The order unit is $\bold 1
=[R_R]$. An {\it order-ideal} in $\Cal V(R)$ is a submonoid $S$ of $\Cal V(R)$
that is order-hereditary (i\.e\. if $x\leq y$ and $y\in S$, then $x\in
S$). We say that two order-ideals $S$ and $T$ of $\Cal V(R)$ are {\it
  orthogonal} provided that $S\cap T=0$. We denote by $\Cal L(R)$ the
lattice of (two-sided) ideals of $R$, and by $\Cal L(\Cal V(R))$ the lattice of
order-ideals of $\Cal V(R)$. If $I$ is a (proper) two-sided ideal of $R$,
we denote by $\Cal F\Cal P(I,R)$ the class of modules $P$ in $\Cal F\Cal P(R)$
such that $PI=P$, and by $\Cal V(I)$ the set of isomorphism classes of
elements from $\Cal F\Cal P(I,R)$. It can be seen that $\Cal V(I)$ depends
only on the structure of $I$ as a ring without unit. (In fact, we can also
describe $\Cal V(I)$ as the monoid of equivalence classes of idempotents in
$\Bbb M_{\infty}(I)$.) By construction, $\Cal V(I)$ is an order-ideal of
$\Cal V(R)$. If $I$ and $J$ are ideals of $R$, then it is clear that
$\Cal V(I\cap J)=\Cal V(IJ)=\Cal V(JI)=\Cal V(I)\cap \Cal V(J)$. In
particular, if $I$ and $J$ are orthogonal ideals of $R$ we observe that $\Cal
V(I)$ and $\Cal V(J)$ are orthogonal order-ideals of $\Cal V(R)$.

If $R$ is an exchange ring, then the map $\phi:\Cal L(R)\to \Cal L(\Cal V(R))$ given
by $I\to \Cal V(I)$ is a surjective lattice homomorphism. Indeed, if $S$ is
an order-ideal of $\Cal V(R)$, then if we denote by $I(S)$ the ideal of $R$
generated by the set $\{e=e^2\mid [eR]\in S\}$, we have that
$\phi(I(S))=\Cal V(I(S))=S$. This correspondence is an isomorphism if we
restrict the domain to the lattice $\Cal L_\Cal J(R)$ of {\it semi-primitive}
ideals, i\.e\. those ideals $I$ of $R$ such that $\Cal J(R/I)=0$, which
form a lattice with infima given by intersections and suprema by
Jacobson radicals of sums (see \cite{{\bf 22}}).
\enddefinition

\medskip

\proclaim{8.7. Theorem} Let $R$ be a unital ring, and consider the
three conditions: \roster
\item "(i)" $R$ is a $QB-$ring\;;
\item "(ii)" If $n\cdot\bold 1+b_1=n\cdot \bold 1+b_2$ in $\Cal V(R)$ for
  some $n\ge 1$, then  we can find
  orthogonal order-ideals $S_1$ and $S_2$ in $\Cal V(R)$ and elements $x, c_1,
  c_2$, such that $c_i\in S_i$ for $i=1,2$, and moreover
  $b_1+c_1=b_2+c_2$ and $x+c_1=n\cdot\bold 1 = x+c_2$\;;
\item "(iii)" If $a+b_1=a+b_2$ in $\Cal V(R)$ then there exist orthogonal
  order-ideals $S_1$ and $S_2$ in $\Cal V(R)$ and elements $c_1, c_2$, such
  that $c_i\in S_i$ for $i=1,2$, and  $b_1+c_1=b_2+c_2$.
\endroster
We always have (i)$\implies$(ii)$\implies $(iii), but if $R$ is a
semi-primitive exchange ring all three conditions are equivalent.
\endproclaim

\demo{Proof} (i)$\implies $(ii). Choose representatives $B_1$ and
$B_2$ for $b_1$ and $b_2$ such that $A_1\oplus B_1 = A_2\oplus B_2$
for some $A_i\simeq R_R^n$. Applying Condition (ii) in Lemma 8.5 we
find decompositions $A_i=N\oplus C_i$ and $B_1\oplus C_1\simeq
B_2\oplus C_2$ for some finitely generated projective modules $N,
C_1, C_2$, such that $C_iI_i=C_i$ for a pair of orthogonal ideals
$I_1, I_2$ in $R$. Let $S_i=\Cal V(I_i)$ and put $x=[N]$ and $c_i=[C_i]$.
(Note that $S_1$ and $S_2$ are orthogonal order-ideals of $\Cal V(R)$.)
Then $b_1+c_1=b_2+c_2$ and $x+c_1=n\cdot\bold 1=x+c_2$, and evidently
$c_i\in S_i$.

(ii)$\implies $(iii). This follows from the fact that $\bold 1$ is an
order-unit in $\Cal V(R)$.

Assume now that $R$ is a semi-primitive exchange ring satisfying
(iii). Consider an arbitrary von Neumann regular element $a$ in $R$
with partial inverse $b$, and put $p=ab$ and $q =ba$, so that $aR=pR$
and $Ra=Rq$. Then 
$$ 
\aligned R\oplus (1-p)R &= qR\oplus (1-q)R\oplus
(1-p)R \\ &\simeq pR\oplus (1-p)R\oplus(1-q)R = R\oplus (1-q)R\;.
\endaligned
$$ 
By (iii) there exist $c_1, c_2$ in $\Cal V(R)$, such that
$$
[(1-p)R]+c_1=[(1-q)R]+c_2\;. 
$$ 
Moreover, $c_i\in S_i$, where $S_1$
and $S_2$ are orthogonal order-ideals in $\Cal V(R)$.

Since $R$ is an exchange ring $\Cal V(R)$ is a refinement monoid, cf\.
\cite{{\bf 4}, Corollary 1.3}, so we can find decompositions $$
\aligned &[(1-p)R]=m_{11}+m_{12}\,,\quad c_1=m_{21}+m_{22}\,,\\
&[(1-q)R]=m_{11}+m_{21}\,,\quad c_2=m_{12}+m_{22}\;.
\endaligned
$$

Since $S_1\cap S_2=0$, we have $m_{22}=0$. Corresponding to this we
have pairs of orthogonal idempotents $e_1, f_1$ and $e_2, f_2$ in
$R$, such that $1-p=e_1+f_1\,,\; 1-q=e_2+f_2\,,\; e_1R \simeq e_2R$
and $[f_1R]= c_2$ and $[f_2R]=c_1$. For $i=1,2$ let $I_i$ be the
ideals of $R$ generated by the idempotents $e$ in $R$ such that $[eR]\in
S_i $, respectively. Then $I_1\cap I_2$ contains no non-zero
idempotents, and so $I_1\cap I_2\subset \Cal J(R)$ because $R$ is an
exchange ring. Since $R$ is semi-primitive $I_1\cap I_2=0$.

Choose $c$ and $z$ in $R$ such that $e_1=cz$ and $e_2=zc$. Then $c$
is a regular element in $(1-p)R(1-q)$, so $a\le a+c$ by Lemma 2.7.
Moreover, $$ 1-(a+c)(b+z)=1-ab-cz=1-p-e_1=f_1\in I_2\;, $$ and
similarly $1-(b+z)(a+c)=f_2\in I_1$. Since $I_1\cap I_2=0$ it follows
that $a+c\in R^{-1}_q$.

We have shown that every regular element in $R$ extends to an element
in $R^{-1}_q$, whence $R$ is a $QB-$ring by Theorem 8.4.
\hfill$\square$
\enddemo

\medskip 
 
\example{8.8. Examples} Theorem 8.7 provides us with more interesting
examples of $QB-$rings. 

\noindent{\bf (A)} Let $\Bbb F$ be a field and define $\Bbb M_\infty(\Bbb
F)=\varinjlim \Bbb M_n(\Bbb F)$ as in 6.9 and 7.16. This is a well-known
example of a simple (non-unital) $B-$ring. Now let $\Bbb B(\Bbb F)$ denote the
algebra of all row- and column-finite matrices over $\Bbb F$. This algebra
contains $\Bbb M_\infty(\Bbb F)$ as its only non-trivial ideal. In fact, if
$V=\Bbb F^{(\omega)}$ denotes the countably infinite dimensional vector space
over $\Bbb F$, and $\Bbb B(\Bbb F)$ is regarded as a subalgebra of
$\End_\Bbb F (V)$, then $\Bbb B(\Bbb F)$ is the idealizer of $\Bbb
M_\infty(\Bbb F)$ in $\End_\Bbb F(V)$.

It is known that $\Bbb B(\Bbb F)$ is not a regular ring, but  O'Meara has
proved to us that $\Bbb B(\Bbb F)$ is an exchange ring. (Private communication.)

Since idempotents in $\Bbb M_\infty(\Bbb F)$ are equivalent if and only if
they have the same (finite) rank, and moreover all idempotents in 
$\Bbb B(\Bbb F)\setminus \Bbb M_\infty(\Bbb F)$ are equivalent to $1$, we see
that $\Cal V(\Bbb B(\Bbb F))=\{0,1, \cdots , \infty\}$, where $\infty =
[1]$. Then it is easy to check that Condition (ii) in Theorem 8.7 is
satisfied, and hence $\Bbb B(\Bbb F)$ is a $QB-$ring.

\noindent{\bf (B)} Replacing the field with a $QB-$ring $A$, whose structure
is not too intricate, say $A$ being a simple exchange ring with prescribed
$\Cal V(A)$, it seems safe to predict that $\Bbb B(A)$ will often be a
$QB-$ring. We already know from Corollary 6.10 that $\Bbb M_\infty (A)$ is a
$QB-$ring, and since this is an ideal in $\Bbb B(A)$ we just have to control
the quotient. However, the ideal structure
of $\Bbb B(A)/\Bbb M_\infty (A)$ can be quite complicated (see
\cite{{\bf 3}}, \cite{{\bf 19}}, \cite{{\bf 25}}, \cite{{\bf 26}}).

\noindent{\bf (C)} Following Chen, \cite{{\bf 14}}, we say that a unital ring
$R$ satisfies {\it related comparability} if whenever we have two idempotents
$p$ and $q$ such that $1-p\sim 1-q$, then there is a central idempotent $z$ in
$R$, such that $zpR$ is isomorphic to a direct summand in $zqR$ and $(1-z)qR$
is isomorphic to a direct summand in $(1-z)pR$. It follows from
\cite{{\bf 14}, Theorem 2} that if $R$ has related comparability, and
if $a$ is any von Neumann regular element in $R$, then we can find a
partial inverse $v$ for $a$, and  a central idempotent $z$ such that
$zv$ is left invertible and $(1-z)v$ is right invertible in
$R$. Evidently $v\in R^{-1}_q$. Thus, if $R$ is a semi-primitive
exchange ring with related comparability, then $R$ is a $QB-$ring by
Theorem 8.4.

\noindent{\bf (D)} An exchange ring $R$ is said to satisfy {\it general
comparability} provided that for any modules $A$ and $B$ in $\Cal F\Cal P(R)$, 
there exists a central idempotent $e$ in $R$ such that $Ae$ is
isomorphic to a direct summand of $Be$ and $B(1-e)$ is isomorphic to a 
direct summand of $A(1-e)$. Thus, $R$ satisfies the fundamental
comparison lemma for Murray-von Neumann equivalence of projections in
a von Neumann algebra. Any exchange ring satisfying general 
comparability also satisfies related comparability. Consequently,
every regular, self-injective ring is a $QB-$ring. This class contains
$\End (V_D)$, where $D$ is a division ring and $V$ is a right
$D-$vector space. More generally, it contains $\End (A_R)$, where $A$
is any non-singular, quasi-injective right module over a ring $R$,
cf\. \cite{{\bf 16}, Corollary 1.23}.

\noindent{\bf (E)} We finally consider the rings constructed in \cite{{\bf 5},
\S 3}. Let $R(p,q)$ be the ring construction in Section 3 of \cite{\bf
5} based on a simple unit-regular ring $L$ and idempotents $p$ and $q$
in $L$. Then $R(p,q)$ is a
prime regular ring with a unique maximal ideal $N(p,q)$. It is proved
in \cite{{\bf 5}, Lemma 3.1} that $R(p,q)$ has stable rank one if and
only if $p\sim q$. By using the same techniques and the fact that
$R(p,q)_q^{-1}=R(p,q)_r^{-1}\cup R(p,q)_l^{-1}$ (because $R(p,q)$ is a
prime ring), one can show that $R(p,q)$ is a $QB-$ring if and only if
either $p\lesssim q$ or $q\lesssim p$. In particular, Example 3.2 (in
\cite{{\bf 5}}) provides an example of a stably finite, regular ring
$U$ which is not a $QB-$ring, but satisfies $2-$comparability (see
\cite{{\bf 16}, p. 275}). In fact, $U$ satisfies the stronger
property of almost comparability. This contrasts heavily with
the situation of regular rings satisfying the comparability axiom (in
the sense of \cite{{\bf 16}}), which are all $QB-$rings since they
satisfy general comparability.
\endexample

\example{8.9. Example} Let $S$ be the ring constructed in Example 7.16. Recall
that $S$ is an extension of the two $B-$rings $\Bbb M_\infty(\Bbb F)$ and
$\Bbb F(t)$, the field of rational functions on a countable field 
$\Bbb F$. Denote by $\rho:S\to \Bbb F(t)$ the quotient map. Selecting
appropriate subrings of $\Bbb F(t)$ with stable rank one it is possible to
produce other examples of $QB-$rings which are neither $B-$rings nor exchange
rings. For example, consider the subring $\Bbb F(t)_{1,2}$ of $\Bbb F(t)$ of
those rational functions $f(t)/g(t)$ such that $g(1)$ is non-zero and $g(2)$
is non-zero (assuming that $2\ne 0$ in $\Bbb F$). The subring $S_{1,2}$ of $S$
which is the inverse image of $\Bbb F(t)_{1,2}$ through $\rho$ will be a
$QB-$ring, by the same argument as in Example 7.16. Since $\Bbb F(t)_{1,2}$ is
semilocal but not local we get that $S_{1,2}$ is not an exchange ring. Again
$S_{1,2}$ is an extension of two $B-$rings, $\Bbb M_\infty(\Bbb F)$ and $\Bbb
F(t)_{1,2}$, and $\bsr(S_{1,2})=2$. We can think of $S$ as a nice localization
of $S_{1,2}$.
\endexample

\vskip1truecm
%\vskip0.5truecm
%\newpage

\subhead\nofrills{\bigbf 9. Extremally Rich C*--Algebras}\endsubhead 
\bigskip

Let $A$ be a $C^*-$algebra, i\.e\. an algebra of bounded operators on a
complex Hilbert space, which is closed under taking adjoints and under norm
convergence. Then $A$ is a semi-primitive ring, and the $^*-$operation makes
it particularly simple to link the left and the right structure of $A$.

The concept of quasi-invertibility for elements in a unital $C^*-$algebra was
explored in \cite{{\bf 8}, Theorem 1.1}, and it was shown that $A^{-1}_q$ is the
open subset of $A$ consisting of elements of the form $A^{-1}\Cal E
(A)A^{-1}$. Here $\Cal E (A)$ denotes the set of extreme points in the unit
ball of $A$, identified by Kadison as the partial isometries $v$ in $A$ such
that 
$$
(1-vv^*) \perp (1-v^*v)\,,
$$
cf\.\, \cite{{\bf 23}, 1.4.7}.

In \cite{{\bf 8}, \S 3} a unital $C^*-$algebra $A$ was defined to be 
{\it extremally rich} if the set of quasi-invertible elements in $A$ was dense in
$A$, i\.e\. if $(A^{-1}_q)^= = A$. This seemed a most appropriate
generalization of Rieffel's notion of $C^*-$algebras with topological stable
rank one, since by results of Rieffel \cite{\bf 28}, and of Herman and
Vaserstein \cite{\bf 18} a $C^*-$algebra has Bass stable rank one if and only
if $(A^{-1})^= = A$, i\.e\. $A$ has topological stable rank one.

To show that our notion of $QB-$ring reduces to extremal richness in the case of
$C^*-$alge\-bras we offer the following result, of which the second half
essentially is a modification of the argument from \cite{\bf 18}.

\medskip

\proclaim{9.1. Proposition} For an element $x$ in a unital $C^*-$algebra $A$
the following conditions are equivalent:
\roster
\item "(i)" $x \in (A^{-1}_q)^= $\;;
\item "(ii)" $x \in \cl (A^{-1}_q)$\;.

\endroster
\endproclaim

\demo{Proof} (i) $\implies$ (ii) If $ax+b=1$ for some $x, b$ in $A$, then 
$$
1 = (ax+b)^*(ax+b) \le 2 x^*a^*ax + 2 b^*b \le 2\Vert a\Vert^2 x^*x + 2b^*b\;.
$$
Thus, $x^*x+b^*b \in A^{-1}$. By Theorem 3.3 in \cite{\bf 8} there is an
extreme point $u$ in $\Cal E(A)$ such that $x+ub\in A^{-1}_q$. In particular,
$x\in \cl (A^{-1}_q)$.

\noindent (ii)$\implies$ (i)  If $x\in \cl (A^{-1}_q)$ and $\varepsilon > 0$ is
given, define $b=(1-\varepsilon ^{-2}x^*x)_+ $, the positive part of the
element $1-\varepsilon^{-2}x^*x$. Then by spectral theory
$$
x^*x+b = x^*x + (1-\varepsilon^{-2}x^*x)_+  \ge \varepsilon^2 1\;,
$$
since $t+(1-\varepsilon^{-2}t)_+ \ge \varepsilon^2$ for every real $t\ge
0$. Consequently $x^*x+b\in A^{-1}$. By assumption we therefore have
$x+yb \in A^{-1}_q$ for some $y$ in $A$. Take a natural number $n$ and define
$$
z=(x+yb)(1+nb)^{-1} \quad\text{in}\quad A^{-1}_q\;.
$$
Then with $|x|=(x^*x)^{\tfrac 12}$ we have
$$
\aligned
\Vert x-z\Vert & \le \Vert x(1-(1+nb)^{-1})\Vert + \Vert y\Vert
\Vert b(1+nb)^{-1}\Vert\\
& = \Vert |x|\,nb(1+nb)^{-1}) \Vert +  
\tfrac 1n \Vert y\Vert \Vert nb(1+nb)^{-1}\Vert\;.
\endaligned
$$
Regarding $b$ as a function of $|x|$ it vanishes on sp$(|x|)$ if
$t\ge\varepsilon$, so $\Vert |x|\,b\,\Vert \le \varepsilon$. Moreover, 
$\Vert \,nb(1+nb)^{-1}\Vert \le 1$ by spectral theory, whence 
$$
\Vert x-z\Vert \le \varepsilon +  \tfrac 1n \Vert y\Vert < 2\varepsilon
$$
for $n$ large enough. Thus, $x\in (A^{-1}_q)^=$. \hfill$\square$
\enddemo
 
\proclaim{9.2. Corollary} In a $C^*-$algebra $A$ we always have\newline
\newline

\qquad\qquad\qquad $\cl (A^{-1}_q) = \text{\rm cr} (A^{-1}_q) \quad (=
(A^{-1}_q)^=)$\;. \hfill$\square$\endproclaim

\medskip

A (unital) $C^*-$algebra $A$ has {\it real rank zero} provided that every
self-adjoint element can be approximated arbitrarily well by
self-adjoint, invertible elements. This concept was introduced and
explored in \cite{{\bf 7}}, where a number of equivalent
characterizations were also provided. 

To appreciate the following result, which was originally proved by L.G. Brown
and the second author by a direct argument [unpublished], note that if $v$ and
$w$ are partial isometries such that $w$ {\it extends} $v$ in the sense of
2.4, then we actually have $v=vv^*w=wv^*v$, which are the relations usually
employed to describe the relation $v\leq w$. To see this, assume that
$v=vxv=vxw=wxv$ for some $x$. Then $v^*v\le w^*w$, and therefore
$$
\aligned
v^*v =v^*x^*w^*wxv &=v^*x^*(v^*v+(1-v^*v))w^*w((1-v^*v)+v^*v)xv \\
&=(v^*v+v^*x^*(1-v^*v))w^*w((1-v^*v)xv+v^*v)\\
&=v^*v+v^*x^*(1-v^*v)w^*w(1-v^*v)xv\;.
\endaligned
$$ 
It follows that $w(1-v^*v)xv=0$, i\.e\. $v=wxv=wv^*vxv=wv^*v$.

\proclaim{9.3. Proposition} A unital $C^*-$algebra $A$ of real rank zero is
extremally rich if and only if every partial isometry in $A$ extends to an
extreme partial isometry.
\endproclaim

\demo{Proof} The necessity of the condition is proved in \cite{{\bf 8},
  Proposition 2.6}, cf\.  \cite{{\bf 8}, Corollary 4.3}, so it only remains to show
sufficiency. 

Towards this end let $a$ be a von Neumann regular element in $A$. Then $a$ has
a polar decomposition $a=v|a|$, and the spectrum of $|a|$ has a gap $]0,
\varepsilon[$ for some $\varepsilon > 0$. Therefore $e=1-v^*v+|a|$ is
an invertible, positive element in $A$, and $x=e^{-1}v^*$ is a partial inverse
for $a$, since $axa=v|a|e^{-1}v^*v|a|=vv^*v|a|=v|a|=a$. 

By assumption $v$ extends to an extreme partial isometry $u$ in $A$, so that 
$(1-uu^*)A(1-u^*u)=0$, and $v=uv^*v=vv^*u$. Set $b=ue$ and
$y=e^{-1}u^*$. Then $1-by=1-uu^*$ and $1-yb=e^{-1}(1-u^*u)e$, so $b\in
A^{-1}_q$ with quasi-inverse $y$. Moreover,
$$
\aligned
&axb=v|a|e^{-1}v^*ue=vv^*ue=ve=v|a|=a\,,\\
&bxa=uee^{-1}v^*v|a|=uv^*v|a|=v|a|=a\,,
\endaligned
$$
so $a\le b$. Since $a$ was arbitrary, it follows from Theorem 8.4 and
\cite{{\bf 4}, Theorem 7.2} that $A$ is a $QB-$ring, and therefore
extremally rich by Proposition 9.1. \hfill$\square$
\enddemo

\example{9.4. Remarks} Using examples from the theory of extremally rich
$C^*-$algebras we can exhibit several phenomena in the theory of
$QB-$rings. Thus, it is possible to have a unital ring $R$, with an idempotent
$p$, such that both $pRp$ and $(1-p)R(1-p)$ are $QB-$rings, but $R$ is not,
cf\. \cite{{\bf 8}, Example 6.12}. 

In the example above the ideals $I_1$ and $I_2$ of $R$, generated by $p$ and
$1-p$, respectively, will still be $QB-$rings, and of course
$I_1+I_2=R$. Since these are complex algebras we also have that $\widetilde I_1 =
I_1+\Bbb C 1$ is a $QB-$ring. This shows that we can not in Corollary 7.15
relax the condition on the subring from being a $B-$ring to being a
$QB-$ring. 

In \cite{{\bf 8}, Example 1.3} we find a $C^*-$algebra $A$ which is not a
$QB-$ring although it is a direct limit of algebras isomorphic to $\Bbb B(\Cal
H)$. However, the embedding of one copy of $\Bbb B(\Cal H)$ into the next is
not prime. This shows that our condition for prime embeddings in Proposition
5.7 is necessary. 

The algebra $A$ contains an element $a$ such that $\pi(a)$
is left or right invertible in $\pi(A)$ for each primitive (i\.e\. 
irreducible) representation $(\pi, \Cal H)$ of $A$. However, $a\notin
\cl(A^{-1}_q)$ by \cite{{\bf 8}, Corollary 1.10}. We can therefore not in
general hope to simplify the definition of quasi-invertibility to just a
question of one-sided invertibility in prime or primitive quotients.
\endexample

%\newpage

\enddocument